\documentclass[11pt]{amsart}

\usepackage{amssymb}
\usepackage{latexsym}
\usepackage{amsmath}
\usepackage[dvips]{graphics}
\usepackage{graphicx}
\usepackage{verbatim}
\usepackage{psfrag}

\newenvironment{dem}{\textbf{Proof.}\par}
{\begin{flushright}$\Box$\end{flushright}}

\newtheorem{definition}{Definition}[section]
\newtheorem{proposition}{Proposition}[section]
\newtheorem{lemma}{Lemma}[section]
\newtheorem{theorem}{Theorem}[section]
\newtheorem{corollary}{Corollary}[section]

\def\Pa{{\mathcal P}}

\def\Ta{{\mathcal T}}
\def\nn{{\noindent}}

\def\bO{{\bf O}}
\def\ot{\otimes}

\def\Z{{\mathbb Z}}
\def\C{{\mathbb C}}
\def\Tr{\mbox{Tr}}
\def\Ua{\mathcal U}
\def\ra{\rightarrow}
\def\Ma{\mathcal  M}

\def\Ha{\mathcal  H}
\def\N{{\mathbb N}}
\def\R{{\mathbb R}}
\def\ts{\times}
\def\cxm{\C\langle (U_i,U_i^*, A_i)_{1\le i\le m}\rangle}
\def\cxA{\C\langle ( A_i)_{1\le i\le m}\rangle}

\def\D{\Delta}
\def\bt{{\bf t}}
\def\bk{{\bf k}}
\def\bl{{\bf \ell}}
\def\e{\varepsilon}
\def\un{1}
\def\sur{\genfrac{}{}{0pt}{}}
\def\ppq{\leqslant}

\def\bU{{\bf U}}
\def\mun{\hat{\mu}^N}
\def\munv{\mu^N_V}

\newcommand{\MNC}{\mathcal{M}_N}

\newcommand{\UNC}{\mathcal{U}_N}
\newcommand{\ONC}{\mathcal{O}_N}
\def\ONCm{\ONC^m}
\def\UNCm{\UNC^m}
\newcommand{\ANC}{\mathcal{A}_N}
\newcommand{\ANR}{\mathcal{A}_N^1}

\def\l{\lambda}
\def\UU{U}
\def\MM{\mathbb{M}}
\def\TM{{\widetilde {\mathbb{M}}}}
\def\L{\mathbb{L}}
\def\muf{\mu^f}
\def\Pa{\mathcal P}

\begin{document}

\title{Asymptotics of unitary and orthogonal  matrix integrals}

\author{Beno\^\i{}t Collins}
\address{CNRS, ICJ,
Lyon 1 Claude Bernard University, 43 bd du 11 Novembre 1918, 69622
Villeurbanne, France
and Dept. Math. Stat., University of Ottawa, collins@math.univ-lyon1.fr}
\author{Alice Guionnet}
\address{UMPA, CNRS  UMR 5669, ENS Lyon,46 all\'ee d'Italie,
69007 Lyon, France, aguionne@umpa.ens-lyon.fr}
\author{Edouard Maurel-Segala}
\address{LPMA, Universit\'es Paris VI, 175 rue du Chevaleret,
75013, Paris, emaurel@proba.jussieu.fr}


\begin{abstract}
In this paper,  we prove that in small parameter regions,
 arbitrary unitary matrix integrals converge in the large $N$ limit
and match their formal expansion.
Secondly we give a combinatorial model for our matrix integral
asymptotics
and investigate examples related to free probability and the HCIZ integral.
Our convergence result also leads us to new results
 of smoothness of microstates. We finally generalize our approach
to integrals over the othogonal group.
\end{abstract}

\maketitle

\section*{Introduction}

Matrix integrals provide   models for  physical systems
(2D quantum gravitation, gauge theory, renormalization, etc...),
and
generating series  for a wide family
of  combinatorial objects (see e.g \cite{Ho,Zv}).

Gaussian integrals are the most studied. It was shown
by Br\'ezin,
Itzykson, Parisi and Zuber \cite{BIPZ}
that perturbations of Gaussian  integrals
expand formally as a generating function of maps,
sorted by their genus when the dimension $N$ of the matrices
is regarded as a parameter. Such `topological' expansions were
shown also to hold in the large $N$ limit,
and then to
  match with the formal expansion  on a
mathematical level of
rigor by two authors
\cite{GMa1,GMa2,Ma}
and previously in the one matrix case
in  \cite{APS,ACKM}
and \cite{EM}. The relation of Gaussian
matrices with the enumeration of maps
is
  an easy consequence of Wick calculus  -or equivalently
Feynman diagrams- see \cite{Zv} for a
 good introduction.
 According to 't Hooft \cite{Ho}, such
topological expansion should hold in the more general context of
models invariant under
unitary conjugation. This leads us to concentrate in this
article  on
 matrix integrals  given by
\begin{equation}
\label{unitarymatrixint}
I_N(V, A_i^N):=\int_{\UNCm} e^{N \Tr (V(U_i,U_i^*,A_i^N, 1\le i\le m))}
 dU_1\cdots dU_m\end{equation}
where
$(A_i^N, 1\le i\le m)$ are $N\ts N$ deterministic  uniformly
bounded  matrices, $dU$ denotes the
Haar
measure on the unitary group $\UNC$ (normalized so that $\int_{\UNC}dU=1$)
and
$V$ is a  polynomial function
in the  non-commutative variables $(U_i,U_i^*,A_i^N, 1\le i\le
m)$. $\Tr$ denotes the usual trace on $N\ts N$ matrices given by
$\Tr(A)=\sum_{i=1}^N A_{ii}$.

We will study  in this article
the first order asymptotics of   matrix integrals given by
\eqref{unitarymatrixint} when  the joint distribution of the
$(A_i^N, 1\le i\le m)$
converges; namely for all polynomial
function $P$ in $m$ non-commutative indeterminates
\begin{equation}
\label{limitA}
\lim_{N\ra\infty}
\frac{1}{N}\Tr(P(A_i^N, 1\le i\le m))=\tau(P)
\end{equation}
for some linear functional $\tau$ on
the set of polynomials.
Without loss of generality, we will assume
that
$(A_i^N, 1\le i\le m)$ are Hermitian matrices.

For technical reasons, it is  convenient to assume
  that
the polynomial $V$
is such
that $\Tr (V(U_i,U_i^*,A_i^N, 1\le i\le m))$ is real for
all $U_i\in \UNC$, all Hermitian matrices $A_i^N$,  for all
$i\in\{1,\cdots,m\}$ and
$N\in\N$ .

Under those very general assumptions, the only result
proved so far is  the formal  convergence
of these matrix integrals. Namely, it was proved in \cite{c1}
 by one author that
for each $k$, the quantity
$$\frac{\partial ^k}{\partial z^k}N^{-2}\log \int_{\UNCm} e^{zN Tr
 (V(U_i,U_i^*,A_i^N, 1\le i\le m))} dU_1\cdots dU_m|_{z=0}$$
converges towards a constant $f_k (V,\tau)$ depending only on the limiting
distribution of the $A_i^N$'s and $V$. Besides, if $V$ is polynomial with integer coefficients,
$f_k (V,\tau)$ is a polynomial function with integer coefficients of the limit moments of the $A_i^N$'s.

\goodbreak
In this paper we will answer affirmatively to the following, previously open questions:
\begin{enumerate}
\item
Does the limit of the matrix integrals exist for small parameters $z$?
\item
Does  the power series $\sum_k \frac{z^k}{k!}f_k(V,\tau)$ have a strictly positive radius of convergence?
\item
Is the limit of the matrix integral equal to the power series?
\end{enumerate}

The following Theorem is a precise decription of our results:
\begin{theorem}
\label{main}
Under the above hypotheses
and if we further assume that
the spectral radius of the
matrices
$(A_i^N, 1\le i\le m, N\in\N)$
is uniformly bounded (by say $M$),
there exists $\varepsilon=\varepsilon(M,V)>0$
so that for $z\in [-\varepsilon, \varepsilon]$, the limit
$$F_{V,\tau}(z):=\lim_{N\ra\infty}\frac{1}{N^2}
\log \int_{\UNCm} e^{z N \Tr (V(U_i,U_i^*,A_i^N, 1\le i\le m))}
 dU_1\cdots dU_m$$
exists.
Moreover, $F_{V,\tau}(z)$ is an analytic
function of $z\in \C\cap B(0,\varepsilon)=\{z\in\C:|z|\le
\varepsilon\}$
and for all $k\in\N$,
$$\left.\frac{\partial ^k}{\partial z^k} F_{V,\tau}(z)\right|_{z=0}
= f_k (V,\tau).$$
\end{theorem}
This also implies that
the series
 $F_{V,\tau}(z)$ has a positive radius
of convergence, a result
which had not been proved by
the techniques of \cite{c1}
based on Weingarten
functions.

Our approach is based
on non-commutative
differential
calculus (in particular on the resulting Schwinger-Dyson or Master loop
equations) and perturbation  analysis as
developed in the context of Gaussian
matrices in \cite{GMa1,GMa2,Ma}. Another possibility
to prove the equality between real and formal limits would
have been
 to show convergence
of the integrals for complex parameters $z$.
We have not yet been able to
follow  this line successfully, and this remains an open
question.

An important example of unitary matrix integral is the so-called
spherical integral, studied by Harish-Chandra and
by Itzykson and Zuber,
$$HCIZ(A,B):=\int_{U\in\UNC}e^{N \Tr (U^*AUB)}dU.$$
This integral is of fundamental importance in analytic Lie theory
and was computed for the first time by Harish-Chandra in \cite{HC}.
In the last two decades it has also become an issue to study
its large dimension asymptotics \cite{GZ,ZZ,cs,GM1}.

Theorem  \ref{main}
holds true for the HCIZ integral.
It thus relates the results of \cite{c1}
(which  computed the formal limit of  the HCIZ integral)
and those of \cite{GZ}
(where the limit of $HCIZ(A,B)$
was obtained (regardless of any small
parameters assumptions) by using large deviations
techniques). Let us recall the limit found in \cite{GZ}.
Let us define
$$I(\mu)=\frac{1}{2}\mu(x^2) +\frac{1}{2}\int\int\log|x-y| d\mu(x)d\mu(y).
$$
If  $\mu_A$ (resp. $\mu_B$)
denote the limiting spectral measure
of $A$ (resp. $B$), assume that $I(\mu_A)$ and $I(\mu_B)$
are finite. Then,
the limit of $N^{-2}\log HCIZ(A,B)$ is given,
according to \cite{GZ}, by
\begin{align}
I(\mu_A,\mu_B)&:=\lim_{N\ra\infty}\frac{1}{N^2}\log  HCIZ(A,B)
\label{formulaHCIZ}\\
&
=-I(\mu_A)-I(\mu_B)-\frac{1}{2}\inf_{\rho,m}\left\{
\int_0^1 \int\left( \frac{m_t(x)^2}{\rho_t(x)}
+\frac{\pi^2}{3}
\rho_t(x)^3\right)dxdt\right\}
\nonumber
\end{align}
where the infimum is taken over $m,\rho$ so that the measure-valued
process  $\mu_t(dx)=
\rho_t(x)dx$
is a continuous process, $\mu_0=\mu_A$,
$\mu_1=\mu_B$ and  $$\partial_t
\rho_t(x)+
\partial_x m_t(x)=0.$$
The inf over $(\rho_t,m_t)$ is taken  (see \cite{G1})
at the solution of
an Euler equation for isentropic flow
 with negative pressure $-\frac{\pi^2}{3} \rho^3$.

Theorem  \ref{main}
 shows that $I(\mu_{\sqrt{\beta} A},\mu_{\sqrt{\beta}B}
)$ depends analytically on $\beta$
in a real neightbourhood of $0$, a result
which is not obvious from formula \eqref{formulaHCIZ}.
Moreover, the coefficients of this expansion
count certain planar graphs (see section \ref{seccombinatorics}),
as summarized in the following theorem.
\begin{theorem} Denote $\sqrt{\beta}\sharp\mu$ the probability measure
$$\sqrt{\beta}\sharp\mu(f)=\int f(\sqrt{\beta}x) d\mu(x).$$
Assume that $\mu_A$ and $\mu_B$ are two compactly supported
probability measures. Then, there exists $\beta_0>0$
such that for all $\beta\in [0,\beta_0]$,
$$I(\sqrt{\beta}\sharp\mu_A,\sqrt{\beta}\sharp\mu_B)=\sum_{n\ge 0}
\beta^n {\mathbb M}_n(\mu_A,\mu_B)$$
converges absolutely. Moreover, we have
$${\mathbb M}_n(\mu_A,\mu_B)=\sum_{m\mbox{ admissible maps of }\Sigma_n}
M_m(\mu_A,\mu_B).$$
$\Sigma_n$ is the set of planar maps drawn above $n$ vertices defined
as stars of
type $U^*AUB$ by gluing pairwise oriented arrows and
possibly rings
 and $M_m(\mu_A,\mu_B)$ is the weight of the map.

\end{theorem}
We refer the reader to section \ref{seccombinatorics} for
the  definitions of stars, admissible maps and weights.
Our definition of planar maps is more complicated than those arising
in the topological expansion of Gaussian matrix models (and which
are directly related with Wick Gaussian calculus and Feynman diagrams): the sums are signed
and we have a notion of admissibility. However it was
an
open question in mathematical
physics to have a graphical model for
unitary integrals (see \cite{ZZ}). Moreover,
this graphical interpretation
gives a new understanding of
cumulants formulae (see section \ref{applicfreep}).

The convergence of other integrals was still unknown and it is one
of the points of this paper to show their convergence.
We use it to study Voiculescu's microstates
entropy evaluated at a set of  laws which are small perturbations
of the law of free variables, and prove regularity of microstates
\begin{theorem}
For tracial states $\mu$ satisfying
suitable assumptions described in Theorem \ref{thmfreeentropy},
\begin{align*}\chi(\mu)&:=\liminf_{\sur{\e\downarrow 0}{
k\uparrow \infty}}\liminf_{N\ra\infty}
\frac{1}{N^2}\log \mu_N^{\ot m}\left(\Gamma_R(\mu,\e,k)\right)
\\
&=\limsup_{\sur{\e\downarrow 0}{
k\uparrow \infty}}\limsup_{N\ra\infty}
\frac{1}{N^2}\log \mu_N^{\ot m}\left(\Gamma_R(\mu,\e,k)\right)
\end{align*}
and a formula for $\chi(\mu)$ can be given.
\end{theorem}

 This result generalizes section 4 in \cite{GMa1}.

The paper is organized as follows:
after setting our working framework (section \ref{notations}),
we study the action of perturbations
upon the integral $I_N(V,A_i^N)$
and deduce some properties
of the related Gibbs measure; namely that
the so-called empirical distribution
of the matrices under this Gibbs
measure satisfies asymptotically
an equation called the Schwinger-Dyson equation (section \ref{matrixmodels}).
Then, we study this equation
and obtain uniqueness for parameters of the potential $V$ small enough
(section \ref{studySD}) and analyticity (section \ref{formalanalyticity}).

We  also describe a (new) combinatorial solution of Schwinger-Dyson equation
(section \ref{seccombinatorics}) and therefore of the first order
of unitary matrix integrals. We deduce applications
of these results to
free probability (section \ref{secfreeproba}) and
to the convergence of matrix integrals $I_N(V,A_i^N)$ (section \ref{secmatrixintegrals}).
Finally, we point out
some consequence
of our result for free
entropy (section \ref{secfreeentropy}). At last, we consider the case
where integration holds over the orthogonal group instead of the unitary
group and show that the first order of such integrals is the same.

\section{Notations}\label{notations}
Let ${\Ua}_N$ be  the set of $N\ts N$ unitary matrices,
$\Ma_N$ the set of $N\ts N$ matrices
with complex entries,
$\mathcal{H}_N$ the subset of Hermitian matrices of $\Ma_N$ and
$\mathcal{A}_N$ the subset of antihermitian matrices
of $\Ma_N$.
We let $m$ be a fixed integer number
throughout this
article.
We denote by $(A^N_i)_{1\le i\le m}$
a $m$-tuple of $N\times N$ Hermitian matrices.
We shall assume that the sequence $(A^N_i)_{1\le i\le m}$
is uniformly bounded for the operator
norm, and without loss of generality
that they are bounded by one,
$$\sup_{N,i}\|A^N_i\|_\infty=\sup_{N,i}\lim_{p\ra\infty}
\left(\Tr( ( A_i^N)^{2p})\right)^{\frac{1}{2p}}\le 1.$$

\subsection{Free $*$-algebra}
Let $\cxm$ denotes the set of polynomial
functions in the non-commutative indeterminates
$(U_i,U_i^*, A_i)_{1\le i\le m}$ with the relation $$U_iU_i^*=U_i^*U_i=1.$$
Note that in general we may want to consider models with a number
of ``deterministic"  indeterminates $A_i$ different from
the number of
``random unitary" indeterminates $U_i$, but this general case
can be obtained from the previous one by looking only at a sub-algebra and
our convention shortens a little the notations.
The algebra $\cxm$ is equipped with the involution $*$ so that
$ A_i^*=A_i$, $(U_i)^*=U_i^*$; $(U_i^*)^*=U_i$ and for any $X_1,\cdots,X_n\in
(U_i,U_i^*, A_i)_{1\le i\le m}$, any $z\in\C$,
$$(zX_1X_2\cdots X_{n-1}X_n)^*={\bar z}X_n^*X_{n-1}^* \cdots X_2^* X_1^*.$$
Note that for any $U_i\in\UNC$,
$A_i\in\Ha_N$, and $P\in\cxm$,
$$\left( P( U_i, U_i^*, A_i,1\le i\le m)\right)^*
=P^*( U_i, U_i^*, A_i,1\le i\le m)$$
where in the left hand side $*$ denotes the
standard involution on $\Ma_N$.
We denote $\cxm_{sa}$ the set of self-adjoint  polynomials;
$P=P^*$, and  $\cxm_a$ the set of
anti-self-adjoint  polynomials ; $P^*=-P$.
In the sequel, except when something different is explicitly assumed,
we shall make the hypothesis  that the
potential $V$ belongs to $\cxm_{sa}$,
which insures that
$\Tr\left(V((U_i, U_i^*, A_i^N)_{1\le i\le m})\right)
$ is real-valued  for all $U_i\in \UNC$
and $A_i^N\in\Ha_N$. Conversely, any potential
$V$ such that $\Tr\left(V((U_i, U_i^*, A_i^N))\right)
$ is real-valued for all $U_i\in \UNC$
and $A_i^N\in\Ha_N$ is self-adjoint up to the addition of some commutators
(which does not change the trace). Indeed, this implies that $\Tr
\left((V-V^*)((U_i, U_i^*, A_i^N)_{1\le i\le m})\right)$ vanishes, which
insures that $V-V^*=\sum_l P_lQ_l-Q_lP_l$ for some polynomials $P_l,Q_l$, cf \cite{CD07} Lemma 2.9 for a probabilistic
proof or \cite{KS07}, Proposition 2.3 for a direct proof (in the real symmetric case, but
directly adaptable to the Hermitian case).
Then, $W:=V+\sum_l(Q_lP_l-P_lQ_l)/2$ is self-adjoint.

\subsection{Non-commutative derivatives}

On $\cxm$, we define
the non-commutative derivatives $\partial_i$, $1\le i\le m$,
given by the linear form such that
$$\partial_i A_j=0,\quad
\partial_i U_j =1_{i=j} U_j\otimes 1\quad
\partial_i U_j^* =-1_{i=j}  1\otimes  U_j^*,\quad \forall j,$$
and satisfying the Leibnitz rule, that  for $P,Q\in \cxm$,
\begin{equation}\label{leibniz}
\partial_i(PQ)=\partial_i P\times (1\otimes Q)+
(P\otimes 1) \times \partial_i Q.\end{equation}
Here, $\times$ denotes the product $P_1\otimes Q_1\times P_2\otimes Q_2=
P_1P_2\otimes Q_1Q_2$.
We also let $D_i$ be the corresponding {\it cyclic}
derivatives such that if $m(A\otimes B)=BA$,
$D_i=m\circ \partial_i$.

If $q$ is a monomial
in $\cxm$, we more specifically have
\begin{eqnarray}
\partial_i q&=&\sum_{q=q_1U_iq_2}
q_1 U_i\otimes q_2
- \sum_{q=q_1U_i^* q_2}
q_1 \otimes U_i^* q_2\label{partialq}\\
D_i q&=&\sum_{q=q_1U_iq_2}
q_2q_1 U_i
- \sum_{q=q_1U_i^* q_2}
U_i^* q_2q_1.\label{Dq} \\
\nonumber
\end{eqnarray}

\subsection{Bounded tracial states}
Let $\Ta$ be the set of tracial states on
the algebra generated by the variables $ (U_i,U_i^*, A_i)_{1\le i\le m}$,
i.e. the set of
linear
forms on $\cxm$ such that for all $P,Q\in\cxm$,
$$\mu(PP^*)\ge 0,\quad \mu(PQ)=\mu(QP),\quad \mu(1)=1.$$
Throughout this article, we restrict
ourselves to tracial states $\mu\in\Ta$  such that
$$\mu((A_i^*A_i)^n)\le 1\quad\forall n\in\N, \, \forall
i\in\{1,\cdots,m\}.$$
We denote $\Ma$ this subset of $\Ta$.

Note that for any monomial
$q\in\cxm$,
the H\"older's inequality implies that for any $\mu\in \Ma$,
\begin{equation}\label{holder}
\mu(qq^*)\le 1.
\end{equation}
We endow $\Ma$ with its weak topology:
$\mu_n$ converges to $\mu$ if and only if for
all $P\in\cxm$,
$$\lim_{n\ra\infty}
\mu_n(P)= \mu(P).$$
By equation \eqref{holder} and since the above topology is the
product topology,
$\Ma$ is a compact metric space by Banach Alaoglu's theorem.

We denote $\mun
$ the empirical distribution
of matrices $A_i^N\in \Ha_N$ and $U_i\in \Ua_N$
which is given for all $P\in\cxm$ by
$$\mun(P)=\frac{1}{N}\Tr
\left(P( U_i,U_i^*, A_i^N , {1\le i\le m})\right).$$

This object will be of crucial interest for us.

The notation
 $\Ma|_{(A_i)_{1\le i\le m}}$  stands for the set of
tracial states of $\Ma$
restricted to the algebra generated by the $(A_i)_{1\le i\le m}$.
In particular, the limiting distribution $\tau$
given by \eqref{limitA}
belongs to $ \Ma|_{(A_i)_{1\le i\le m}}$.

\subsection{Tracial power states}
Let $V\in\cxm_{sa}$ and  $\mu^N_V$ be the distribution on $\UNCm$
given by
$$\mu^N_V(dU_1,\cdots, dU_m)
=I_N(V,A_i^N)^{-1}
\exp(N\Tr(V))dU_1\cdots dU_m.$$
We define, for all $P\in\cxm$,
$$\bar\mu^N_V(P)
:=E_{\munv}[\mun(P)]:=\frac{\int \frac{1}{N}\Tr P e^{N\Tr V}dU_1\dots dU_n}{\int e^{N\Tr V}dU_1\dots dU_n}.$$
In the following, an $n$-tuple of
monomials $(q_i)_{1\le i\le n}$ in
$\cxm$ will be fixed and we shall take
$V=V_{\bf t}=\sum_{i=1}^n t_i q_i$.
Then, $\bar \mu^N_{V_{\bf t}}(P)$
can be seen as a power series in the $t_i$'s;
\begin{equation}\label{etatf}\bar \mu^N_{V_{\bf t}}(P):=
\sum_{\bk\in
\N^n}\frac{\bt^\bk}{\bk!}\left.\frac{\partial^{|\bk|}}{\prod_i\partial
t_i^{k_i}}\right|_{t_i=0}\frac{E[\mun(P)e^{ N^2
\mun ( V_{\bt})}]}{E[e^{
N^2\mun(
 V_{\bt})}]}.
 \end{equation}
We will call $\mu$ a `tracial power state' of $\Ma$ if and only if it is a map
$$\mu: \cxm\to \C[[\bt]]$$
with for all $a,b$, $\mu (ab)=\mu(ba)$.
Here $\C[[\bt]]$ is the algebra of power series in
the variables $t_1,\cdots,t_n$.
In particular, we may view $\mu^N_{V_{\bf t}}$
as a tracial power state of  $\Ma$.

\subsection{Cumulants.}
The classical cumulants $\{C_k\}_{k\ge 0}$
 are defined via their formal generating function:
$$\log E(e^{tX})=\sum_{k\geq 0} t^kC_k(X,\ldots ,X)/k!$$
This equality holds also in a complex neighborhood of $0$ for $t$ if $X$ is bounded.
We also define the cumulants $C_\bk$ for $\bk$ in $\N^n$:
$$\log E(e^{t_1X_1+\cdots+t_nX_n})=\sum_{\bk\in \N^n} \bt^\bk C_\bk(X_1,\ldots ,X_n)/\bk!$$
where $\bk=(k_1,\cdots,k_n)$, $\bk!=\prod_i k_i!$, $|\bk|=\sum_i k_i$ and $\bt^\bk=\prod_i t_i^{k_i}$.
Note that:
$$C_\bk(X_1,\ldots ,X_k)=C_{|\bk|}(X_1,\cdots,X_1,\cdots,X_n,\cdots,X_n)$$
where in the previous list the variable $X_i$ appears
$k_i$ times.

Let us recall some properties of these cumulants.
\begin{proposition}\label{technical}
The following two statements hold true:
\begin{enumerate}
\item
$$\frac{E(Ye^{t_1X_1+\ldots +t_nX_n})}{E(e^{t_1X_1+\ldots +t_nX_n})}=\sum_{\bk\in \N^n} \bt^\bk C_{1,\bk}(Y,X_1,\ldots ,X_n)/\bk!$$
\item
\begin{align*}
&\frac{E(YZe^{t_1X_1+\ldots +t_nX_n})}{E(e^{t_1X_1+\ldots +t_nX_n})}-\frac{E(Ye^{t_1X_1+\ldots +t_nX_n})}{E(e^{t_1X_1+\ldots +t_nX_n})}\frac{E(Ze^{t_1X_1+\ldots +t_nX_n})}{E(e^{t_1X_1+\ldots +t_nX_n})}\\
&=
\sum_{k\geq 0} \bt^\bk C_{1,1,\bk}(Y,Z,X_1,\ldots ,X_n)/{\bf k}!
\end{align*}
\end{enumerate}
\end{proposition}

\begin{dem}
Item (1) is obtained by replacing $t_1X_1+\ldots +t_nX_n$ by $yY+t_1X_1+\ldots +t_nX_n$ and differentiating
the generating function of the cumulants
in $y$ at $y=0$.

Item (2) is obtained by replacing $tX$ by $yY+zZ+tX$ and differentiating
the equality defining the cumulants in $y$ and $z$  at $y,z=0$.
\end{dem}

\section{Matrix models}\label{matrixmodels}

We first investigate the asymptotic behavior of the
random state $\mun$  under
$\mu^N_V$ as a random tracial state.
We then consider $\bar\mu^N_V=\mu^N_V(\mun)$
evaluated at a polynomial and study its convergence as
a formal power series in the parameters of the potential $V$.
We show that they satisfy asymptotically the same type
of equations called Schwinger-Dyson (or Master loop)
equations.

\subsection{Behavior of $\mun$}

The main result of this section
is the following
\begin{theorem}
\label{mainsec1}
Assume that $V$ is self-adjoint.
For all polynomial \\$P\in\cxm$,
$$\lim_{N\ra\infty}\left\lbrace
\mun\otimes\mun (\partial_i P)
+\mun(D_i VP)\right\rbrace=0\quad \munv \,  a.s.$$
In particular, any limit point $\mu\in\Ma$
of $\mun$ under $\munv$  satisfies the
Schwinger-Dyson equation
\begin{equation}
\label{SD0}
\mu\otimes \mu(\partial_i P)+\mu(D_i V P)=0
\end{equation}
for all $P\in\cxm$
and $\mu|_{(A_i)_{1\le i\le m}}=\tau.$
\end{theorem}
The idea of the proof, rather common in quantum field
theory and successfully used in \cite{GMa1,GMa2,Ma},
is to obtain equations on $\mun$ by
 performing an infinitesimal change of
variables in $I_N(V,A_i^N)$.
More precisely we make the
change of variables $\bU=(U_1,\cdots,U_m)\in \UNCm
\to \Psi(\bU)=(\Psi_1(\bU),\cdots,\Psi_m(\bU))\in \UNCm$ with
 $$\Psi_j(\bU)= U_je^{\frac{\l}{N} P_j(\bU)}$$
where the $P_j$ are antisymmetric polynomials (i.e. $P_j^*=-P_j$).
This change of variables  becomes very close to the identity
as $N$ goes to infinity, reason
 why it is called ``infinitesimal''.

\begin{lemma}\label{lempsi}
The function $\Psi$ is a local
diffeomorphism
and its Jacobian $ J_{\Psi}$ has the following expansion
when $N$ goes to infinity
$$|\det J_{\Psi}(\bU)|=e^{\frac{\l}{N} \sum_i\Tr\otimes \Tr(\partial_i
P_i
(U_i,U_i^*,A_i, 1\le i\le m))
+O(1)}$$
where $O(1)$ is uniform on the unitary group (but may depend on $P$).
\end{lemma}

\begin{dem}
Let us first recall the following two elementary results of differential
geometry:
\begin{enumerate}
\item
The map $\exp:\MNC\longrightarrow\MNC$ is differentiable and:
$$\mbox{Diff}_M \exp . H :=\lim_{\varepsilon\ra 0}
\varepsilon^{-1} (e^{ M+\varepsilon H}-e^M)
= \left(\sum_{k=0}^{+\infty}\frac{(\mbox{Ad}_M)^{k}}{(k+1)!} H\right) e^{M}$$
where $\mbox{Ad}_M$ is the operator defined by $\mbox{Ad}_M H=MH-HM$.
\item
If $P\in\cxm$ is considered as a function of the $U_i$'s,
then it is differentiable and its differential with respect to the $i$-th variable
in the direction $A$, for $A$ in $\ANC$,
is
$$\mbox{Diff}_iP. A := \lim_{\varepsilon\ra 0}
\varepsilon^{-1} (P(U_1,\cdots, U_{i-1}, U_ie^{\varepsilon A},U_{i+1}
,\cdots)-P(\bU)) =\partial_i P \sharp A.$$
\end{enumerate}
As a consequence, if we
fix $A$  in $\ANC$ and $i\in\{1,\cdots,m\}$, one has
\begin{align*}
\mbox{Diff}_i \Psi_j(\bU).A&
=1_{i=j}U_jA+U_j \mbox{Diff}_{\frac{\l}{N}P_j(U)}\exp.(\frac{\l}{N}\partial_i
P_j\sharp A)\\
&=1_{i=j}U_jA+\frac{\l}{N}\sum_{k=0}^{+\infty}
U_j\frac{(\mbox{Ad}_{\frac{\l}{N}P_j(U)})^{k}}{(k+1)!}
(\partial_i P_j\sharp A)e^{\frac{\l}{N}
P_j(U)}\\
&=1_{i=j}U_jA+U_j\frac{\l}{N}\Phi_{ij}(\bU) A.
\end{align*}
with $\Phi_{ij}(\bU)$ the linear map from $\ANC$ into $\MNC$ given by
$$\Phi_{ij}(\bU) A:=
\sum_{k=0}^{+\infty}
\frac{(\mbox{Ad}_{\frac{\l}{N}P_j(\bU)})^{k}}{(k+1)!}
(\partial_i P_j\sharp A)e^{\frac{\l}{N}
P_j(U)}.$$
We can factorize the term $U_j$ to obtain
\begin{equation}
\mbox{Diff} \Psi(\bU)=\UU\circ(\mbox{Id}_{\ANC^m}+\frac{\l}{N}\Phi(\bU))
\end{equation}
with $\UU\circ(M_1,\cdots,M_m)=(U_1M_1,\cdots,U_mM_m)$ and $\Phi$  the linear operator from $\ANC^m$ to $\MNC^m$ whose blocks are the $\Phi_{ij}(\bU)$.

Since the operator norms of the $A_i$'s and the $U_i$'s are uniformly bounded in $N$,
the operator norm of $\mbox{Ad}_{\frac{\l}{N}P_j(\bU)}$ as an operator on
$(\MNC,\|.\|_\infty)$ is also bounded. Thus, $\Phi_{ij}(\bU)$ is a uniformly bounded operator from $\ANC$ to $\MNC$.
Thus,  the norm of
$\frac{\l}{N}\Phi(\bU)$ is less than $1/2$ for $N$ large enough. For those $N$,
$\Psi$ is a local diffeomorphism with positive eigenvalues.

We can now compute the  Jacobian of $\Psi$
$$|\det J_{\Psi}(\bU)|:=|\det \mbox{Diff} \Psi(\bU)|=|\det
\UU||\det(I+\frac{\l}{N}
\Phi(\bU))|.$$
It can be easily checked that $|\det \UU|=1$.

Besides, the positivity of the eigenvalues of $I+\lambda \Phi(\bU)/N$ allows us to replace the determinant by the exponential of a trace:
\begin{align*}
|\det J_{\Psi}(\bU)|&=\exp(\Tr\log(I+\frac{\l}{N}\Phi(\bU)))
=\exp\left(-\sum_{p\geq 1}\frac{(-\l)^p}{pN^p}\Tr(\Phi(\bU)^p)\right).
\end{align*}
Note that since $\Phi$ is a bounded operator on $\ANC$,
which is a space of dimension $N^2$,
the $p$-th term in the previous sum is at most of order
$N^{2-p}$. We only look at the terms up to the order $O(N)$.
A quick computation shows that if
$$\varphi:
\begin{array}{ccc}
\ANC&\to&\ANC\\
X&\to&\sum_lA_l X B_l
\end{array}
$$
is considered as a real endomorphism, $\Tr\varphi=\sum_l \Tr A_l\Tr B_l$.
Indeed, if we consider $E({kl}), 1\le k, l\le N$ the canonical basis
of $\ANC$, $$E({kl})_{rj}:=\sqrt{-1}\frac{1_{r=k,j=l}+1_{r=l, j=k}}{\sqrt{2(1+1_{k=l})}}$$
for $k\le l$ and
$$E({kl})_{rj}:=\frac{1_{r=k,j=l}-1_{r=l, j=k}}{\sqrt{2}}$$ for $k\ge l$,
$\Tr\varphi=\sum_{k, l} \Tr( E(kl)^*\varphi(E(kl)))=\sum_l \Tr A_l\Tr B_l$.
This is sufficient to obtain the first term of the Jacobian:

$$
\frac{\l}{N}\Tr(\Phi(\bU))=\frac{\l}{N}\sum_i\Tr(\Phi_{ii}(\bU))
=\frac{\l}{N}\sum_i\Tr\otimes\Tr(\partial_i P_i(U_j,U_j^*,A_j))+O(1)
$$
with $O(1)$ is uniformly bounded on $\UNCm$. Here we used that the
operator norm of $\mbox{Ad}_{\frac{\l}{N}P_j(\bU)}$ is uniformly small.
\end{dem}

Before making the change of variables we show that
$\Psi$ is  a bijection.
\begin{lemma}
For $N$ large enough, $\Psi$ is a diffeomorphism of $\UNCm$.
\end{lemma}

\begin{dem}
First observe that since $\Psi$ is a local diffeomorphism,
its image is open in $\UNCm$.
Besides, since $\UNCm$ is compact and $\Psi$ is continuous, the image
is compact and therefore closed. Thus by connectedness
 of $\UNCm$, and since $\Psi(\UNCm)$
is closed, open and non-empty, $\Psi$ is surjective.

The only property we still need to prove is the injectivity of $\Psi$.
If $\Psi(U)=\Psi(V)$ then for all $j\in\{1,\cdots,m\}$,
$$U_j^*V_j-I=e^{\frac{\l}{N}P_j(U)}e^{-\frac{\l}{N}P_j(V)}-I.$$
Thus, if $N$ is sufficiently large so that $\frac{\l}{N}P_j(U)$ is in a domain where the
function $\exp$ is
$2$-Lipschitz,  we obtain
\begin{align*}
\|U_j-V_j\|_\infty&=\|U_jV_j^*-1\|_\infty=\|e^{\frac{\l}{N}P_j(U)}e^{\frac{-\l}{N}P_j(V)}-1\|_\infty\\
&=\|e^{\frac{\l}{N}P_j(U)}-e^{\frac{\l}{N}P_j(V)}\|_\infty
\le \frac{2|\l|}{N} \|P_j(U)-P_j(V)\|_\infty
\end{align*}
 with $\|.\|_\infty$  the operator norm.
Since $(P_j, 1\le j\le m)$ are uniformly Lipschitz on
$\UNCm$, we conclude that  $\sum_{j=1}^m\|U_j-V_j\|_\infty
$ vanishes for sufficiently large $N$.
\end{dem}

We can now prove Theorem \ref{mainsec1}.

\begin{dem}
Let us define
$$Y^N(P)=\sum_i\left(\frac{1}{N}\Tr(D_iV P_i)+\frac{1}{N}\Tr\otimes\frac{1}{N}\Tr(\partial_i P_i)\right).$$
We expand
 $\Tr V(\Psi(\bU)_i,\Psi(\bU)_i^*, A_i, 1\le i\le m )$ as
\begin{eqnarray}
&&\Tr(V(\Psi(\bU)_i,\Psi(\bU)_i^*, A_i, 1\le i\le m ))
-\Tr(V(U_i,U_i^{*}, A_i, 1\le i\le m))\nonumber\\
&&
=\frac{\l}{N}\sum_j\Tr (D_jVP_j(U_i,U_i^{*}, A_i, 1\le i\le m) )
+O(N^{-1})\label{kol}\\
\nonumber
\end{eqnarray}
and perform the change of
variables ${\bf U}\ra \Psi({\bf U})$
in $I_N(V, A_i^N)$;
\begin{align*}
I_N(V, A_i^N)&:=\int e^{ N\Tr(V(U_i,U_i^{*}, A_i, 1\le i\le m))} dU_1\cdots dU_m\\
&=
\int e^{ N \Tr(V(\Psi(\bU)_i,\Psi(\bU)_i^*, A_i, 1\le i\le m ))}
|\det J_\Psi(\bU)|  dU_1\cdots dU_m\\
&= \int e^{N Y^N(P) +0(1)}  e^{N\Tr(V(U_i,U_i^{*}, A_i, 1\le i\le m)) } dU_1\cdots dU_m\\
\end{align*}
where we used \eqref{kol} and Lemma \ref{lempsi}.
$O(1)$ is of order one independently of $N$
and uniformly on the unitary matrices $(U_1,\cdots, U_m)$.
Thus we have proved that
$$\int e^{N Y^N(P)} d\mu^N_V(\bU)
=O(1).$$
Borel-Cantelli's lemma thus insures that
$$\limsup_{N\ra \infty} Y^N(P)\le 0\quad a.s.$$
and the converse inequality holds by changing $P$ into
$-P$ since $Y^N$ is linear in $P$. This proves the first statement
of Theorem \ref{mainsec1}.
The last result is simply based
on the compactness of $\Ma$ and the fact that
any limit point must then satisfy
the same asymptotic equations
than $\mun$.
\end{dem}

Another consequence of this convergence is
the existence of solutions
to \eqref{SD0} for any self-adjoint potential $V$
(since any limit point of $\mun$ in the compact metric space $\Ma$
will satisfy it) a fact already proved in \cite{Bi}.
Moreover,  since these
solutions are limit points of $\mun$,
they belong to $\Ma$ and in particular $|\mu(q)|\le 1
$ for any monomial $q$.

\subsection{Moments of $\mun$}

In the sequel, we denote by $E$ the expectation with respect to the
Haar measure on the unitary group.
The goal of this section is to show (see Proposition \ref{qwes})
that cumulants also satisfy a
formal version of Schwinger-Dyson equation.
We start with the following lemma:

\begin{lemma}\label{lemma-formal-sd}
One has, for all $i$ all $N$,
all monomials $q_1,\cdots, q_n$
 and all $\bk=(k_1,\cdots,k_n)$ in $\N^n$,

\begin{align*}
&N^2 E\left( \mun\otimes \mun \left(\partial_i P\right)\cdot \left(\mun (q_1)\right)^{k_1}\cdots \left(\mun (q_n)\right)^{k_n}\right)\\
&+\sum_j k_j E\left(\left(\mun( q_1)\right)^{k_1}\cdots \left
(\mun (q_j)\right)^{k_j-1}\cdots\left(\mun (q_n)\right)^{k_n} \mun (D_iq_j\cdot P)\right)=0
\end{align*}

\end{lemma}

\begin{dem}
Following Lemma \ref{lempsi}, we
write down the change of variable
 $$\Psi_i : \bU\to (U_1,\cdots, U_{i-1},
U_ie^{\l P_i(\bU)}, U_{i+1},\cdots,U_m) $$
in the integral
$\int ((\mun q_1)^{k_1}\cdots (\mun q_n)^{k_n})
dU_1\cdots dU_m$, where the integration is taken over the unitary Haar measure.
Its Jacobian satisfies
$$|\det J_{\Psi}(\bU)|=1+\frac{\l }{N}
\Tr\otimes \Tr (\partial_i P)+o(\l ).$$
and we have the expansion
\begin{align*}\Tr (q_j(\Psi(\bU)_i,\Psi(\bU)_i^*,A_i,1\le i\le m ))&=\Tr
(q_j(U_i,U_i^*,A_i,1\le i\le m))\\&\hspace{-43pt}+\l \Tr (D_i q_j\cdot P(U_i,U_i^*,A_i,1\le i\le m))+\l^2o(\l )
\end{align*}
where the $o(\l)$'s are for a given $P$ uniform bounds in $N$.
The first order of the  Taylor expansion of
this change of variables around $\l =0$ proves the claim.
\end{dem}

\begin{proposition}\label{qwes}
As a formal series equality, one has, for all $i$, for all $\bt$,
$$E[\mun\otimes \mun (\partial_i P)e^{ N^2\mun ( V_{\bt})}]
+  E[\mun (D_i V_{\bt}\cdot P)e^{ N^2\mun (V_{\bt})}]
=0.$$
\end{proposition}

\begin{dem}
Multiplying the
equality
of Lemma \ref{lemma-formal-sd}
by $\bt^\bk N^{2|\bk|-2}/\bk!$ and summing over $\bk$ in $\N^n$
gives  the desired identity.
\end{dem}

Finally we  study the large $N$ limit $\muf$  of these formal states (the index $f$ stands for ``formal'').
\begin{theorem}\label{formal-limit}
Let $V_{\bt}$ be the polynomial $\sum_{j=1}^n t_j q_j$.
For all $P$, the sequence $\bar \mu^N_{V_{\bf t}}(P)$ converges as a formal series (i.e. coefficientwise) when $N$ goes to infinity to
some $\muf(P)$.
Besides, $\muf$ satisfies  the family of equations, for all $i$, for all $P$,
$$\muf\otimes \muf (\partial_i P)+ \muf (D_i V_{\bt}\cdot P)=0.$$
\end{theorem}

\begin{dem}
First, we prove the existence of a limit.
By the first item of Proposition \ref{technical}, we can express
$\bar \mu^N_{V_{\bf t}}(P)$ as a sum over cumulants,
$$\bar \mu^N_{V_{\bf t}}(P)=\sum_{\bk\in\N^n}\bt^\bk C_{1,\bk}(\frac{1}{N}\Tr P,N\Tr q_1,\cdots,N \Tr q_n)/\bk!.$$
The limit in $N$, of the $C_{1,\bk}(\frac{1}{N}\Tr P,N\Tr q_1,\cdots,N \Tr q_n)$ was proved to exists in \cite{c1} so that $\muf$ is well defined.

Item (2) from Proposition \ref{technical} implies
\begin{align*}&\frac{E(\frac{1}{N}\Tr P_1 \frac{1}{N}\Tr P_2 e^{N\Tr V})}{E(e^{N \Tr V})}-
\frac{E(\frac{1}{N}\Tr P_1  e^{N\Tr V})}{E(e^{N\Tr V})}
\frac{E(\frac{1}{N}\Tr P_2  e^{N\Tr V})}{E(e^{N\Tr V})}\\
&=\sum_{k\geq 0}\frac{\bt^\bk}{\bk!}
 C_{1,1,\bk}(\frac{1}{N}\Tr P_1,\frac{1}{N}\Tr P_2,N\Tr q_1,\cdots,N\Tr q_n).
\end{align*}
Now, it follows from \cite{c1} that elements on the right hand side have decay $N^{-2}$ so that the coefficientwise limit is zero.
This can be interpreted as a formal convergence of measure result for the states
$\mun$.

The proof of the Theorem follows from this observation and from
Proposition \ref{qwes}.
\end{dem}

\section{Study of the Schwinger-Dyson equation}\label{studySD}
We have shown that the limit points of the matrix model satisfy
the Schwinger-Dyson equation \eqref{SD0}.
The aim of this section is to study this equation and show that it has
a  unique solution.

\begin{definition}\label{defSD}
Let $\tau\in \Ma|_{(A_i)_{1\le i\le m}}$.
A tracial state $\mu\in\Ma$ is said to satisfy Schwinger-Dyson equation
{\bf SD[V,$\tau$]} if and only if for all $P\in \cxA$,
$$\mu(P)=\tau(P)$$
and for all $P\in\cxm$, all $i\in \{1,\cdots,m\}$,
$$\mu\otimes\mu (\partial_i P)+\mu (D_iV \, P)=0.$$
\end{definition}

Let $V\in\cxm$. One can decompose $V$ in a sum
$$V=\sum_{i=1}^n t_i q_i(U_j,U_j^*, A_j,1\le j\le m)$$
with monomial functions $q_i$ and complex numbers $t_i$.
The monomials $(q_i,1\le i\le n)$ will be fixed hereafter.
 We let $D$
be the maximal degree of the monomials $q_i$.

Here we prove that $\tau$ is uniquely defined
provided that the parameters $(t_i,1\le i\le m)$ are small enough.

\begin{theorem}\label{mainsec2}
Let $D$ an integer and $\tau$ a tracial state in $\Ma|_{(A_i)_{1\le i\le m}}$  be given.
There exists $\varepsilon=\varepsilon(D,m)>0$
such that if $|t_i|\le \varepsilon$,
there exists at most one solution
$\mu$ to {\bf SD[V,$\tau$]}.
\end{theorem}
From
this and Theorem \ref{mainsec1}
we deduce the following
\begin{corollary}\label{convcor}
Assume that $V$ is self-adjoint. Let $D$ an integer and $\tau$ a tracial state in $\Ma|_{(A_i)_{1\le i\le m}}$  be given.
There exists $\varepsilon=\varepsilon(D,m)>0$
such that if $|t_i|\le \varepsilon$,
$\mun$ converges almost surely to the unique
solution $\mu$ of the Schwinger-Dyson equation.
Moreover, $\bar\mu^N_V=\mu^N_V(\mun)$ converges as
well to this solution
as $N$ goes to infinity.

\end{corollary}
This result is obvious since Theorems  \ref{mainsec1}
and \ref{mainsec2}
show that $\mun$ has a
unique limit point, and thus converges
almost surely. The convergence
of $\bar\mu^N_V$ is then a direct consequence of
bounded convergence theorem since
$\mun\in\Ma$.

Actually Theorem \ref{mainsec1} and Corollary \ref{convcor}
do not use the assumption that the matrices $(A_i^N, 1\le i\le m)$ are deterministic,
but only  that they are bounded and
have a converging joint distribution. Therefore these
two results extend to the case where these matrices
are random, independent of the $(U_i, 1\le i\le m)$,
and satisfy the above two conditions almost surely.
This observation implies that our result can also encompass the case
of the truncated $GUE$ or other classical bounded matrix models.

We prove now Theorem \ref{mainsec2}.

\begin{dem}
Let $\mu$ be a solution to {\bf SD[V,$\tau$]}.
Note that if we take $q$
a monomial in $\cxm$, either
$q$ does not depend
on $(U_j,U_j^*,1\le j\le m)$ and then $\mu(q)=\tau(q)$
is uniquely defined, or $q$ can be written as
$q=q_1 U_i^{a} q_2$ for some $i\in \{1,\cdots,m\}$,
$a\in \{-1,+1\}$
and monomials $q_1,q_2$ (Here, note that $U_i^{-1}=U_i^*$).
Then, by the traciality assumption, $\mu(q)=\mu( q_2q_1 U_i^{a})
= \mu(  U_i^{a}q')$ with $q'=q_2q_1$. Remark that we can assume without
loss of generality that the last letter of $q'$ is not $U_i^{-a}$.
We next use {\bf SD[V,$\tau$]} to compute
$\mu(  U_i^{a}q)$ for some monomial $q$.
We assume first that $a=-1$.
Then, by \eqref{leibniz},
$$\partial_i \left( U_i^*q\right)=
-1\otimes (U_i^*q)+ U_i^*\otimes 1 \times \partial_i q.$$
Taking the expectation, we thus find
by \eqref{partialq}, since $\mu(1)=1$,  that
\begin{eqnarray}
\mu(U_i^*q)&=&\mu\otimes\mu( U_i^*\otimes 1 \times \partial_i q)
+\mu(D_i V q)\nonumber\\
&=&
 \sum_{ q=q_1 U_iq_2 }
\mu( q_1 )\mu (q_2)
-\sum_{ q=q_1 U_i^*q_2 }
\mu(U_i^* q_1 )\mu (U_i^*q_2)\nonumber\\
&&+\sum_j t_{ij} \mu( q_{ij} q)\label{eqSD}\\
\nonumber
\end{eqnarray}
where $(t_{ij},q_{ij})$ are such that
 $D_i V= \sum_j t_{ij}q_{ij}$.
Note that the sum runs at
most on $Dn$ terms and that
all the $t_{ij}$ are bounded by $\max |t_i|$.
A similar formula is found when
$a=+1$ by differentiating  $qU_i$ (or by using
$\overline{\mu(qU_i))}=\mu((qU_i)^*)=\mu(U_i^* q^*)$).

We next show that \eqref{eqSD} and its equivalent for $a=-1$
characterize uniquely $\mu\in\Ma$
when the $t_{ij}$ are small enough. It will
be crucial here that $\mu(q)$ is bounded
independently of the $t_i$'s (here
by the constant $1$).

Now, let $\mu,\mu'\in\Ma$  be two solutions to {\bf SD[V,$\tau$]}
and set
$$\D(\ell)=\sup_{\mbox{deg}({q})\le \ell}
|\mu(q)-\mu'(q)|$$
where the supremum holds
over monomials of $\cxm$
with total degree in
 the $U_j$ and $U_j^*$ less than $\ell$.
Namely, if the monomial (or word) $q$
contains $a_j^+$ times $U_j$
 and $a_j^-$ times $U_j^*$ , we assume
$\sum_{j=1}^m (a_j^++a_j^-)\le \ell$.
Note that by traciality of $\mu$,
\begin{eqnarray}
\D(\ell)&=&\max_{\sur{1\le i\le m}{ a\in\{+1,-1\}}}\sup_{\mbox{deg}{q}\le \ell-1}
|\mu( U_i^a q)-\mu'(U_i^a q)|\label{ineq1}\\
\nonumber
\end{eqnarray}
and that by \eqref{eqSD},
we find that, for $q$ with degree less than $\ell-1$,
\begin{eqnarray*}
|\mu( U_i^* q)-\mu'(U_i^* q)|
&\le &
\sum_{ q=q_1 U_iq_2 }
|(\mu-\mu')(q_1 )|
+\sum_{ q=q_1 U_iq_2 }
|(\mu-\mu')(q_2)|
\\
&&\hspace{-3cm}+\sum_{ q=q_1 U_i^*q_2 }
|(\mu-\mu')(U_i^* q_1 )|
+\sum_{ q=q_1 U_i^*q_2 }|(
\mu-\mu') (U_i^*q_2)|
\nonumber\\
&&\hspace{-3cm}+\sum_j t_{ij} |(\mu-\mu')( q_{ij} q)|.
\end{eqnarray*}
A similar formula holds for $|\mu( U_i q)-\mu'(U_i q)|$ by
conjugation, and therefore

$$\D(\ell)
\le 2\sum_{p=1}^{\ell -2}  \D(p) + 2\sum_{p=1}^{\ell -1}\D(p)
+nD\varepsilon \D(\ell+D-1)$$
where we used that $\deg(q_1)\in \{0,\cdots,\ell-2\}$,
$\deg(q_2)\in \{0,\cdots,\ell-2\}$ (but $\D(0)=0$)
and
 $\mbox{deg}(q_{ij})\le D$
and assumed $|t_i|\le \varepsilon$.
Hence, we have proved that

$$\D(\ell)\le 4\sum_{p=1}^{\ell -1}\D(p)
+nD\varepsilon \D(\ell+D).$$
Multiplying these inequalities by $\gamma^\ell$
we get, since $H(\gamma):=\sum_{\ell\ge 1} \gamma^\ell \D(\ell)$ is
finite
for $\gamma<1$,
$$H(\gamma)\le \frac{\gamma}{1-\gamma} H(\gamma)
+\frac{nD\varepsilon }{\gamma^{D}} H(\gamma)$$
resulting with $H(\gamma)=0$ for
$\gamma$ so that
$1>\frac{\gamma}{1-\gamma}+\frac{nD\varepsilon }{\gamma^{D}}$.
Such a $\gamma>0$ exists when $\varepsilon$ is small
enough.
This proves the uniqueness.

\end{dem}

As a corollary, we characterize
asymptotic freeness by a Schwinger-Dyson
equation, a result which was already
obtained in \cite{Voi6}, Proposition 5.17.

\begin{corollary}\label{asfree}
A tracial state $\mu$ satisfies {\bf SD[0,$\tau$]} if and only if
, under $\mu$, the algebra generated by $(A_i,1\le i\le
m)$ and $(U_i,U_i^*, 1\le i\le
m)$ are free
and
the $U_i$'s are two by two  free and satisfy
$$\mu( U_i^a)=0\quad \forall a\in\Z\backslash \{0\}.$$
\end{corollary}

\begin{dem}
By the previous theorem,
it is enough to verify that the law $\mu$
of  free variables $(A_i, U_i,U_i^*)_{1\le i\le m}$ satisfies
{\bf SD[0,$\tau$]}. So take $P=U_{i_1}^{a_1}B_1\cdots
U_{i_p}^{a_p}B_p$
with some $B_k$'s  in the algebra generated
by $(A_i,1\le i\le m)$. We wish to show
that for all $i\in\{1,\cdots, m\}$,
$$\mu\otimes \mu( \partial_i P)=0.$$
Note that by linearity, it is enough
to prove this equality when $\mu(B_j)=0$
for all $j$.
Now, by definition, we have
$$\partial_i P=\sum_{k:i_k=i, a_k>0}\sum_{l=1}^{a_k}
U_{i_1}^{a_1}B_1\cdots B_{k-1} U_i^l \otimes U_i^{a_k-l} B_k\cdots
U_{i_p}^{a_p}B_p\qquad$$
$$\qquad -
\sum_{k:i_k=i, a_k<0}\sum_{l=0}^{a_k-1}
U_{i_1}^{a_1}B_1\cdots B_{k-1} U_i^{-l} \otimes U_i^{a_k+l} B_k\cdots
U_{i_p}^{a_p}B_p.$$
Taking the expectation on both sides, since $\mu(U_j^i)=0$ and $\mu(B_j)=0$
for all $i\neq 0$ and $j$, we see that freeness implies
that the right hand side is null (recall here that
in the definition of freeness, two consecutive elements have to
be in free algebras but the first and the
last element can be in the same algebra).
Thus,
$\mu\otimes\mu(\partial_i P)=0$
which proves the claim.
\end{dem}

\section{Formal solution and analyticity}\label{formalanalyticity}
We have shown
in Theorem \ref{formal-limit}
 that the limit points of the formal model also satisfy
an equation similar to Schwinger-Dyson's
 equation. The only difference is that
one of these equations is on the space of tracial states while the other one is on the space of tracial power states.
In order to prove that the formal model matches the matrix model
we need to study this formal equation
and show that the series have a positive radius
of convergence, hence providing a
solution to {\bf SD[V,$\tau$]} as defined in Definition \ref{defSD}.

\begin{definition}
Let $V_{\bt}=\sum_i t_i q_i$ be a polynomial.
Let $\tau$ be a tracial power state in $\Ma|_{(A_i)_{1\le i\le m}}$.
A tracial power state $\mu\in\Ma$ is said to satisfy Schwinger-Dyson equation
{\bf SD$^f$[$V_\bt$,$\tau$]} if and only if for all $P\in \cxA$,
$$\mu(P)=\tau(P)$$
and for all $P\in\cxm$, all $i\in \{1,\cdots,m\}$,
$$\mu\otimes\mu (\partial_i P)+\mu (D_iV_\bt \, P)=0.$$
(Here, both terms of the above equality are elements of
$\C[[\bt]]$
and the equality is formal.)
\end{definition}

We already know, due to Theorem \ref{formal-limit},
 that there exists a solution to this equation. We now prove that this solution is unique.

\begin{theorem}\label{formal-uniqueness}
There exists a unique tracial power state $\bt\ra \mu_{\bt}$ which satisfies
Schwinger-Dyson equation {\bf SD$^f$[$V_\bt$,$\tau$]}.
\end{theorem}

\begin{dem}
Let $\mu_\bt$ be a
 tracial power state solution of {\bf SD$^f$[$V_\bt$,$\tau$]}.
There exists a family $\mu^{\bk},\bk=(k_1,\cdots,k_n)
\in\N^n$ in the algebraic dual of
$\cxm$ such that for all $P$,
$$\mu_{\bt}(P)=\sum_{\bk\in \N^n}\prod_{i=1}^n\frac{t_i^{k_i}}{k_i!}
\mu^{\bk}(P).$$
We will now show that the $\mu_\bk$ are uniquely inductively defined by the relation given by {\bf SD$^f$[$V_\bt$,$\tau$]}. Let us define $1_j$ the vector
in $\N^n$ which vanishes
on every coordinate except the $j$-th which is $1$.
We get the following equalities, for all $\bk$,
\begin{enumerate}
\item
If $P$ is in $\cxA$,  $\mu^\bk(P)=\tau(P)1_{\bk=0}$,
\item
If $P=RU_iS$ with $S$ in $\cxA$, $\mu^\bk(P)=\mu^\bk(SRU_i)$,
\item
If $P=RU_i^*S$ with $R$ in $\cxA$ and $S$ does not contain any $U_j$
(but may contain the $U_j^*$), $\mu^\bk(P)=\mu^\bk(U_i^*SR)$,
\item
If $q$ does not contain any $U_j$,
\begin{align*}
\mu^\bk(U_i^{*}q)&=
-\sum_{ q=q_1 U_i^{*}q_2 }\binom{\bk}{\bk'}\sum_{\bk'+\bk''=\bk}
\mu^{\bk'}(U_i^{*} q_1 )\mu^{\bk''} (U_i^*q_2)\\
&
+\sum_j k_j \mu^{\bk-1_j}( U_i^*q D_iq_{j}).
\end{align*}
\item
And for all $q$,
\begin{align*}
\mu^\bk(qU_i)&=
-\sum_{ q=q_1 U_iq_2 }\sum_{\bk'+\bk''=\bk}\binom{\bk}{\bk'}
\mu^{\bk'}(q_1 U_i)\mu^{\bk''} (q_2U_i)\\
&+\sum_{ q=q_1 U_i^*q_2 }\sum_{\bk'+\bk''=\bk}\binom{\bk}{\bk'}
\mu^{\bk'}(q_1)\mu^{\bk''} (q_2)
-\sum_j k_j \mu^{\bk-1_j}( D_iq_{j} q U_i).
\end{align*}
\end{enumerate}
One can see that this allows to compute uniquely any $\mu^\bk(P)$.
The first relation takes care of the non random case, the relations 2 and 3  use the traciality
to place a variable $U$ in a convenient place. Finally relations 4 and 5 allow to compute
$\mu^\bk(P)$ as a function which depends on the $\mu^{\bk'}(Q)$ with $\deg Q<\deg P$ and $\bk'\ppq\bk$ (first terms) or on the $\mu^{\bk'}(Q)$ with
$\bk'<\bk$ (last term). This is a well founded induction. Thus the
$\mu^\bk$ are uniquely defined.
\end{dem}

We next show that this solution is not only formal but that it
gives a family of solutions $\mu_{\bt}$ of
the non-formal equation ${\bf SD[V_{\bt},\tau]}$, which depends analytically on the
parameters $(t_i)_{1\le i\le n}$.

\begin{theorem}\label{regularity}
There exists $\varepsilon>0$   such that
for  $\bt\in \C^n, \max_{1\le i\le n} |t_i|\le
\e$, the formal solution $\mu_\bt$ of ${\bf SD^f[V_{\bt},\tau]}$ is indeed a convergent series.
For all polynomials $P$, $\bt\in
B(0,\varepsilon)=\{t\in\C^n:\max_{1\le i\le n}|t_i|\le \varepsilon\}
\longrightarrow \mu_{\bt}(P)$
is analytic.

In other
words,  there exists a family ($\mu^{\bk},\bk=(k_1,\cdots,k_n)
\in\N^n$) in the algebraic dual of $\cxm$ such that for all $P$,
$$\mu_{\bt}(P)=\sum_{\bk\in \N^n}\prod_{i=1}^n\frac{t_i^{k_i}}{k_i!}
\mu^{\bk}(P)$$
converges absolutely for $\max_{1\le i\le n} |t_i|\le
\e$.
\end{theorem}
An immediate consequence of this result is to deduce that the formal
solution is a real solution of  ${\bf SD[V_{\bt},\tau]}$
in a small parameters region, and therefore by Theorem \ref{mainsec2},
equals the real solution. This will be a key to prove
Theorem \ref{main} (see section \ref{secmatrixintegrals}).
\begin{corollary}\label{corSDf}
For small $\bt$, the formal solution of Schwinger-Dyson equation ${\bf
SD^f[V_{\bt},\tau]}$
 converges as a series.
In addition, it matches the real solution of ${\bf SD[V_{\bt},\tau]}$
which thus depends analytically in the parameters $\bt$ of the
potential
in a neighborhood of the origin.
\end{corollary}

Let us now prove  Theorem \ref{regularity}.

\begin{dem}
According to the proof of Theorem \ref{formal-uniqueness} the $\mu^\bk$ are uniquely defined
by the family of relations (1)-(5). We only need to control
 the growth of the coefficients $\mu^\bk(P)$ to show that $\mu_\bt(P)$
is indeed convergent for small enough parameters.

To  bound these quantities,  we  use the Catalan  numbers
$$C_0=1, C_{k+1}=\sum_{0\ppq p\ppq k}C_pC_{k-p}$$
and  the fact that they do not explode too fast;
$C_{k+1}\ppq 4C_k$. We denote
 $C_{\bk}:=\prod_i C_{k_i}$ and
$D_{k}:=A^{k-1}C_{k-1}$ for $k\ge 1$,
$D_0:=0$. The two key properties of this sequence
is first that it is sub-geometric ($D_{k+1}\ppq 4AD_k$)
and secondly it satisfies
$D_k=A\sum_{0<p<k}D_pD_{k-p}.$
Now our induction hypothesis is that there exists $A,B>0$
such that for all $\bk$, for all monomial $P$ of degree $p$,
\begin{equation}\label{boundcoef}
\frac{|\mu^{\bk}(P)|}{\bk!}\ppq C_{\bk}B^{\bk}D_p.
\end{equation}
We prove this bound by induction, and the relations (1)-(5) which
define
the $\mu^\bk$.  For $\bk=(0,\cdots, 0)$ this bound is satisfied since
$D_p\ge 1$.
 We will check the induction  for a polynomial of the form
$qU_i$ since it is the most complicated case.
\begin{align*}\frac{|\mu^{\bk}(qU_i)|}{\bk!}
&\ppq\sum_{ \sur{q=q_1 U_iq_2}{\bk'+\bk''=\bk}}
\frac{|\mu^{\bk'}(q_1 U_i)|}{\bk'!}
\frac{|\mu^{\bk''} (q_2U_i)|}{\bk''!}\\
&+\sum_{\sur{q=q_1 U_i^*q_2}{\bk'+\bk''=\bk}}
\frac{|\mu^{\bk'}(q_1)|}{\bk'!}\frac{|\mu^{\bk''} (q_2)|}{\bk''!}
+\sum_{k_j\neq 0} \frac{|\mu^{\bk-{\un_j}}( D_iq_{j} q)|}{(\bk-\un_j)!}
\end{align*}
Now we use the induction hypothesis.
If $q$ is of degree $p-1$,
\begin{align*}\frac{|\mu^{\bk}(qU_i)|}{\bk!C_{\bk}B^{\bk}D_p}
&\ppq2\sum_{ \sur{0<q<p}{\bk'+\bk''=\bk}}
\frac{C_{\bk'}B^{\bk'}D_qC_{\bk''}B^{\bk''}D_{p-q}}
{C_{\bk}B^{\bk}D_p}\\
&+D\sum_j \frac{C_{\bk-\un_j}B^{\bk-1}D_{p+D}}
{C_{\bk}B^{\bk}D_p}\\
&\ppq 2\prod_i \frac{C_{k_i+1}}{C_{k_i}}\frac{1}{A}
+nD \frac{(4A)^D}{B}.
\end{align*}
The point is that we can choose $A,B>0$ such that this last quantity
is lesser than $1$. For example take $A>4^{n+1}$ and then
$B>2nD(4A)^D$.

Thus,
for $\|t\|:=\max_i |t_i|<1/4B$, for all $P$ in $\cxm$,
the series $\sum_{\bk}\prod_i\frac{t_i^{k_i}}{k_i!}\mu^{\bk}(P)$
is absolutely convergent.
\end{dem}

\section{Combinatorics.}\label{seccombinatorics}

The purpose of this section is to provide a graphical approach to the
solution of the Schwinger-Dyson equation, and therefore to the computation
of unitary matrix integrals and free entropy (see sections \ref{secfreeproba}, \ref{secmatrixintegrals}
and \ref{secfreeentropy}).
Actually, the proof of Theorem \ref{formal-uniqueness}
gives a recursive way of computing formal solutions to the Schwinger-Dyson
equation, and therefore numerical solutions with arbitrary precision.

Before giving a detailed description of our combinatorial model,
 we start with an overview.
We need the notions of a {\bf star},
which is a pictorial encoding of a monomial
of $\cxm$, of {\bf  root star}, which is a distinguished
star,  and of a {\bf map}, which is a
specific planar decoration over a set of stars and one  root star.

The goal of this section is to show that  the limits of integrals
 on the space of unitary matrices are generating function
of the number of some maps as described
above.
However we are not interested in all maps, but rather
on some that arise from an admissible
construction, which leads us to the
concept of {\bf admissible maps}.
Last, we
need
the notion of {\bf weight} of a map,
and our result will be in terms of sum
over admissible maps of  weights.

Let us point out that for the sake of clarity, although our
 natural playground is the algebra
$\cxm$ and our definitions work in full generality, we
restrict ourselves in the examples
to the case of one single unitary matrix $U$ and two variables $A_1=:A$ and
$A_2=:B$.
We first start with the definition of  stars and  root stars,
in the spirit of  \cite{GMa1,GMa2}.

\begin{figure}[ht!]
\psfrag{de}{\begin{huge}{distinguished element}\end{huge}}
\psfrag{re}{\begin{huge}{root element}\end{huge}}
\psfrag{a}{\begin{huge}A\end{huge}}
\psfrag{b}{\begin{huge}B\end{huge}}
\psfrag{ho}{\begin{huge}oriented edge\end{huge}}
\psfrag{o}{\begin{huge}orientation\end{huge}}
\psfrag{r}{\begin{huge}ring\end{huge}}
\begin{center}\resizebox{!}{7cm}{\includegraphics{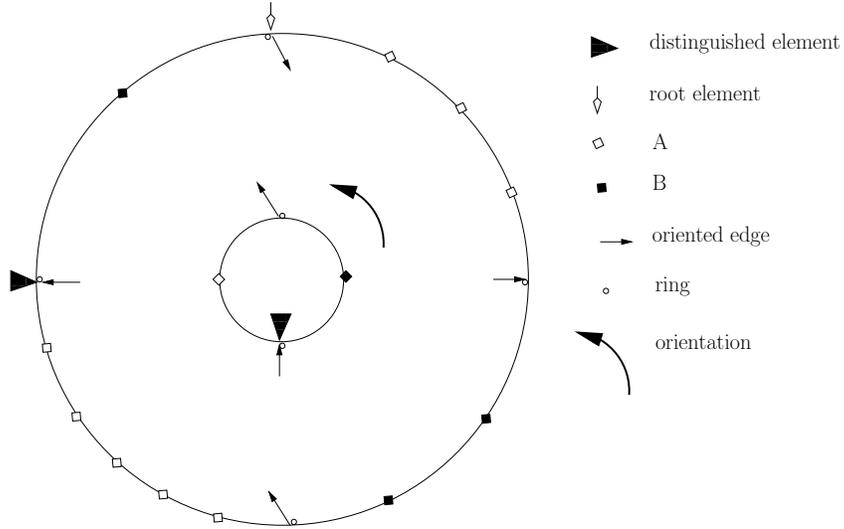}}
\end{center}
\caption{Star of type $U^*AUB$ and root star of type
$U^*A^5UB^2U^*A^3UB$.}\label{vertex}
\end{figure}

\begin{definition}
\begin{enumerate}
\item
A {\bf star }
 is  a circle endowed with the clockwise  orientation,
 decorated with elements such as  colored incoming or outgoing
arrows,  and colored diamonds. One of the element is marked.
\item
To each letter $X_i$
in the alphabet
$(A_i, U_i,U_i^*)_{1\le i\le m}$,  we associate  (bijectively) an
{\bf element}  as follows;
a diamond  of color $i$ if $X_i=A_i$ and
a ring of color $i$ if $X_i=U_i$ or $U_i^*$; in the case of $U_i$
(resp. $U_i^*$)
we attach before the ring an outgoing arrow of  color
$i$
(resp. we attach after the ring an incoming arrow of  color
$i$)
outside of the circle.

\item
To a  monomial  $q\in\cxm$, we associate in a canonical way
a {\bf star  of type $q$} by drawing on the clockwise oriented  circle
the elements associated to  the  successive
letters of $q$, while the element corresponding to the first letter
of $q$ is  marked (or distinguished).

\item
A {\bf root star of type $q$ }  is    obtained
by drawing on the oriented  circle
the elements associated to  the  successive
letters of $q$ in the counter clockwise order, the arrows being drawn
inside
the circle. Its first element is distinguished.
Although
the maps are on the sphere, in the graphical representation
of this section we will draw them on the plane and, to highlight
the role of the root star we will draw it in this section such that it contains all the other stars.
It can be viewed as the star centered in infinity or as the outer face of the dual map.
Besides, on a root star we will distinguish a root element.  If $q$ contains no $U_i$ nor $U_i^*$, there are no root element. If $q$ contains a $U_i$,
the  ring associated to the last $(U_i,1\le i\le m)$
is the root element. If $q$ contains no $U_i$ but some $U_i^*$,
the ring associated to the  first  $(U_i^*,1\le i\le m)$  is called the root element.

\item
A {\bf multistar }
 is a set of $k$   stars  inside
 a root star
drawn on the same plane
with a coherent orientation.

\end{enumerate}

The figure \ref{vertex} shows
a concrete example of a multistar. In the middle of the picture there is a star of type $U^*AUB$ and, surrounding it, a root star of type
$U^*A^5UB^2U^*A^3UB$.
\end{definition}

We are now  ready to introduce the
main objects in our combinatorial model,
namely, maps:

\begin{definition}
A {\bf map} is a decoration of a multistar into
 a connected graph embedded in the plane
by drawing two species of edges between rings:
\begin{enumerate}
\item
A first category of edges, called ``dotted edges", can be drawn
 between rings  either attached to two outgoing
arrows of the same color or to two incoming arrows of the same color.
 These edges can only have rings as end points, not diamonds or arrows.
 Rings can have any number of dotted edges going out of them, possibly none.
\item
A second category of edges, called ``colored oriented edge" arises from the connection of
an arrow going out of a star (associated with a variable $U_i$)
  into  an incoming arrow
   (associated to a variable $U_i^*$) of the same color.
  These colored oriented edges is a pairing between the set of $U_i$'s and the set of $U_i^*$'s: exactly one incoming arrow is glued to each outgoing arrow.
\end{enumerate}
In addition, all the above edges  do not cross, all arrows
are paired but rings can be  attached to any number of
dotted edges (including
to none).

\end{definition}

In the remainder of this section we keep considering pictures drawn on
 the sphere
(and in fact on the plane).
They therefore give rise to graphs with vertices, edges and faces - together
with additional decoration.
For our forthcoming definitions, we need to
clarify the notion of `face': we
 consider that faces of a graph are the connected components
of the complementary of the graph on the sphere.
However, we take the convention that the original stars are
`fattened vertices'. Therefore the interior
of stars will not be considered as faces (neither is the exterior of
the root
star).

Each `face'  component of a map is isomorphic to a
disc;
thus this is an actual face.
This is due to the fact that our map is embedded into a sphere.
This condition would not be granted
 in the case of an embedding into a higher
genus oriented 2D compact manifold.
In this case it would have to stand in the definition of a map of
 `higher genus':
this will be of use for future work but for the sake of
 simplicity we do not emphasize this notion
 in this paper.

Next, we define the weight of a map.
The boundary of a face is homeomorphic to a circle,
it is   given an  orientation (the orientation of the sphere) and
is decorated with diamonds (note that all
arrows have been paired); it thus has the structure
of a star except for the distinguished element.

\begin{definition}
Assume we are given the tracial state $\tau$ of \eqref{limitA}.
\begin{itemize}
\item
First we define the weight of the faces of a map.
The boundary of a face have the structure
of a star, i.e. it has the topology of a circle with some diamonds on it.
We can therefore associated each of these boundaries
with a monomial in the $A_i$'s, given up to cyclic permutation (or
equivalently up to knowing its first letter).
The weight of a face is the trace $\tau(q)$ (which does not depend on
cyclic
permutations)
of the monomial $q$ associated with
its boundary.

 \item
The weight of the map $m$, denoted by $M_{m}(\tau)$, is
the product of the weights
of its faces times a sign given
by $-1$
to the power the number of dotted edges.

\end{itemize}
\end{definition}

As we said before, not all maps will
contribute and we need to
define now the notion
of admissible maps.
Admissibility can be checked
by an inductive procedure
 {\bf IP}, which ressembles Tutte's surgery \cite{Tu}
and which amounts to check one after the other
whether edges of the map
are admissible. Once an edge has been
checked, it is frozen and we continue
by checking the other edges.

{\bf Inductive Procedure {\bf IP} }:

a- If the root star has no root element,
then it can not be connected to any other star
and hence the graph can not be a map unless there is
no other star in which case the map
is just the trivial graph with no edges.

b- The  root star has a root  element which is associated
to a $U_i$ (resp.
a  $U_i^*$), for some $i\in\{1,\cdots,m\}$.

1-Then, we first check the admissibility of the dotted edges
starting from this root element.
We first consider the dotted edge which is the farthest
from  the arrow and declare it admissible
if its other  vertex
  is a ring of an
outgoing  (resp. ingoing) arrow and that  {\it there is no other
dotted edge } attached to this ring which is farther
(amongst the unfrozen dotted edges)
from its arrow. Once this condition is verified, we freeze this
dotted edge and
the  root element
remains
the root element.
We check all dotted edges
of the  root element inductively. Once a dotted edge has been checked
to be admissible, it is frozen and we go on checking the others.
Once  all the  dotted edges of the root element have been checked,
they are frozen and may separate the graph into subgraphs. Thus,
the map may have been cut into
disjoint subgraphs whose
boundary (which may contain dotted edges) is homeomorphic
to a disc (In the case where it has edges glued with an internal
star, we see these other stars as part of the external
star by following all the graph connected to
the external boundary).
 In each of these subgraphs,
we declare the first (following the orientation of the plane) element (corresponding to a $U_i$ or a $U_i^*$)
 after
the last frozen  dotted edge of its boundary
as distinguished. We then define the root element of the boundaries of
these subgraphs by the same procedure as for the root star.  The boundary of
each subgraph  is then a star  and these subgraphs
have now the structure of a map; we will call them submaps.

\begin{figure}[ht!]
\begin{center}\resizebox{!}{4cm}{\includegraphics{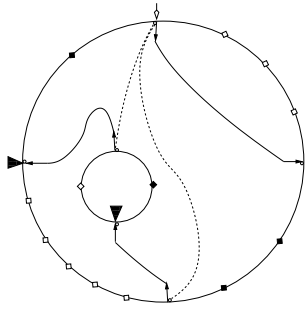}}
\end{center}
\caption{A possible map. Its weight is
$\tau^{\otimes 5}(A^6\otimes B\otimes B^2\otimes A^3 \otimes B)$}\label{map1}
\end{figure}

For instance, in figure \ref{map1}, once the two dotted edges have
been checked
and frozen, the map is cut into two disjoint submaps, the left one
having a fixed frozen dotted edge between the root element and the
inner star (thus forbidding other edges to cross it and allowing us
to consider the inner star as part of the external star). The boundary
of this left subgraph is now seen as a star of type
$q=U^*A^5 UBU^*AUB$.  This left subgraph
has the same distinguished element as before but a new root element (here
the outgoing arrow on its boundary corresponding to the first $U$
in $q$). 

2-
When  all dotted edges are  frozen,
we check that the arrow of the root  element
is paired with an  arrow of the opposite direction (note that
if the root element comes from a $U_i^*$, it can only be paired
with an element of another star since by definition there
is no more outgoing arrows on the root star). The oriented edge
is seen as a fat edge. In particular, if  the oriented edge
link the root star with another star, we see this other
star as part of the root star for the next
step, i.e we identify the root star of type
$QU_iP$ glued to the star of type $RU_i^* S$
(by the marked $U_i$'s) with the root star
of type $PQSR$ with, by convention, the distinguished element
chosen to be the closest element after the glued $U_i^*$.
If the oriented edge link two rings of the root star, two disjoint
 subgraphs are formed and we proceed as in -1-.

c- We continue the inductive procedure on the submaps until all edges
have been checked to be admissible.

Now we can define weighted sum of admissible maps.

\begin{definition}
Assume we are given the tracial state $\tau$ of \eqref{limitA}.

We define the weighted sum of admissible maps constructed
above the stars  $r_1,\cdots,r_n$ and the root star
$P$:
$$\MM_{r_1,\cdots,r_n}(P)=\sum M_{m}(\tau)$$
where the sum runs over all admissible maps $m$
 constructed above $r_1,\cdots,r_n$ with root star $P $.
Assuming that $V_{\bt}=t_1q_1+\ldots +t_nq_n$
 where $q_i$ are monomials,
we define the formal
series:
$$\MM_{\bt}(P)=\sum_{\bk\in\N^n}\frac{\bt^\bk}{\bk!}\MM_{\bk}(P)$$
with $\MM_{k_1,\cdots,k_n}(P)=\MM_{q_1,\cdots,q_1,\cdots,q_n,\cdots,q_n}(P)$
where the monomial $q_j$ appears in $k_j$
successive position
and ${\bf t}^{\bf k}=\prod t_i^{k_i}$, ${\bf k}!=\prod k_i!$.

\end{definition}

Remark that we do not count all the maps which contain the stars
 $r_1$,\dots,$r_n$ but only those that are constructed using our inductive
rules; they
for instance forbid to glue the two
same rings more than
twice.

However, a given map is counted
at most  once since there is only one way to decompose
it using the procedure {\bf IP}.
 Indeed, it is easy to check that at each step we have only one possibility for the next step
since the dotted edges have to be drawn one after
the other following the orientation
and no new dotted edge can be drawn after
the arrow of the root has been glued.

{\bf Example}

\begin{figure}[ht!]
\begin{center}\resizebox{!}{4cm}{\includegraphics{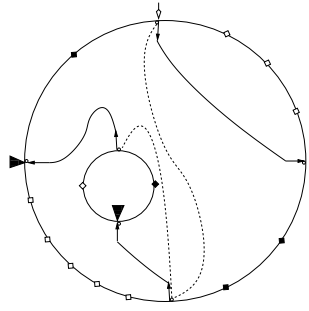}}
\end{center}
\caption{Another one. Its weight is
$\tau^{\otimes 5}(A^6\otimes B\otimes B^2\otimes A^3 \otimes B)$}\label{map2}
\end{figure}

\begin{figure}[ht!]
\begin{center}\resizebox{!}{4cm}{\includegraphics{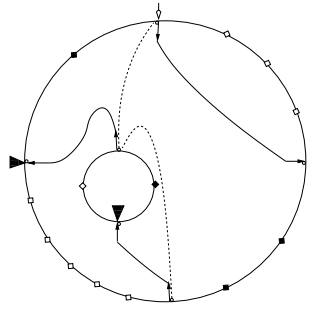}}
\end{center}
\caption{A counterexample: {\bf IP} is violated because of the order of the
dotted edges at the root element}\label{map3}
\end{figure}

\begin{figure}[ht!]
\begin{center}\resizebox{!}{4cm}{\includegraphics{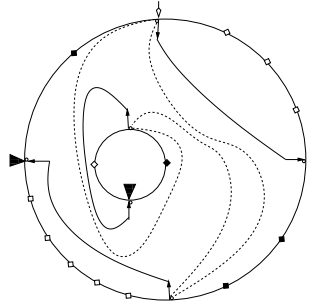}}
\end{center}
\caption{An admissible map. Its weight is $\tau^{\otimes 6}(A^5\otimes
A\otimes  B\otimes B^2\otimes A^3 \otimes B)$
 }\label{map4}
\end{figure}

Let us show some examples. We start from one root
star and a star
 on the sphere (see figure \ref{vertex}).
We want to construct maps above these stars with our rules, starting with the root element
 shown by the arrow outside the root star. Figures \ref{map1}, \ref{map2} and
\ref{map4} are examples of such maps. Note that the weights of the maps of figures
\ref{map1} and \ref{map2} are the same, the only difference is the way
the three rings are glued. There is a third way to glue those three
rings shown in figure \ref{map3} which is a
 map but can not be obtained by our rule of construction (and thus is
not admissible).

We now come to the main theorem of this section, namely
the graphical expansion result for $\MM_{\bf t}$:

\begin{theorem}\label{maincombinatorics}
Let $V=\sum_{1\ppq i\ppq n} t_i q_i$ be a polynomial.
Let $\mu_{\bf t}$ be a solution of ${\bf SD}[V_{\bf t},\tau]$ and $\MM_{\bf t}$ be  the
formal series defined for monomials $P$ by
$$\MM_{\bf t}(P)=\sum_{\bk\in \N^n}\prod_{i=1}^n\frac{t_i^{k_i}}{k_i!}
\MM_{\bk}(P)$$
where $\MM_{\bk}(P)$ is the weighted sum of planar maps with one
root  star of type $P$
and $k_i$ stars  of type $q_i$. If we extend the definition of $\MM_{\bf t}$ by linearity to
any polynomial $P$ then the series $\MM_{\bf t}(P)$ is
absolutely convergent in a neighborhood
of the origin and,
$$\MM_{\bf t}(P)=\mu_{\bf t}(P).$$
\end{theorem}

\begin{dem}

For the sake of clarity we first prove the case $V=0$
and show that $\MM(P):=\MM_{\bf 0}(P)=\mu_{\bf t}(P)$ for a monomial $P$.

We proceed by induction on the total  degree in  $U_i$,
$1\le i\le m$,  in  $q$.

Suppose that there is no variable $U_i$ in $P$.
 Then either there is no variable $U_i^*$
and both sides of the equality are equal to $\tau(P)$, or there
is a $U_i^*$ and both sides vanish: the left hand side by freeness between $U_i$ and the $A_i$'s and the fact that all non-trivial moments of $U_i$ is $0$
and the right hand side because one can not glue the arrow
coming out from this $U^*_i$ anywhere.

We assume our  identification proved  when the degree of $P$ in the $U_i$'s
is less than $k$.
We next take $q$ with degree in the $U_i$'s equal to
$k+1$.
Thus we can assume that there is a $U_i$ in $P$, and we consider
the last one in $P$ so that $P=pU_i b$ with
$b$ a polynomial in the $U_j^*$ and the $A_j$'s,
$1\le j\le m$. By definition, $\MM(pU_ib)=\MM(bpU_i)$
since it depends only on the position of the last $U_i$.
Thus, we may assume  that $P$ is of the form
$QU_i$ with $Q$ of degree
$k$. We apply Schwinger-Dyson equation to this quantity:

\begin{equation}\label{sanspotentiel}
\mu(QU_i)=-\sum_{Q=RU_iS}
\mu(RU_i)\otimes \mu(SU_i)+\sum_{Q=RU_i^*S}
\mu(R)\otimes \mu(S)
\end{equation}

Now, we can apply our induction
hypothesis since all polynomials
appearing in the right hand side
have degree strictly smaller than $k+1$.

We need to show that this is exactly the induction relation for maps.
To construct a map above a star of type $QU_i$, we first look at
the root element $U_i$ and we have to decide
what to do first with the dotted edges.
There are two possibilities:
\begin{enumerate}
\item
The first possibility is that there is no
 dotted edge going outside of the ring of the root.
 In such a case, we can glue
 the arrow to any other arrow  of opposite direction and
of the same color (corresponding to a variable $U_i^*$).
This implies that $Q$ decomposes into $RU_i^*S$ and we construct an
oriented edge between $U_i$ and $U_i^*$. Thus we separate the map into
 two parts and we have to construct a map above the $R$ part and another one above the $S$ part
(this is the case 2 of {\bf IP}).
This gives
$$\MM(R)\MM(S)$$
possibilities which is exactly the possibilities
counted by the second term in the right hand side of \eqref{sanspotentiel}.
\item
The second possibility is that we glue the root ring to another ring with a dotted edge.
Thus $Q$ must decompose into $RU_iS$
and the creation of the  dotted edge  amounts to decompose the
map
into $RU_i$ and $SU_i$ and again to continue the construction of the map
we will have to construct a map above the $RU_i$ part and another one
above the $SU_i$ part (note here that
when a dotted edge is attached to a circle of an $U_i$,
 the arrow and the circle keep their structure and live on the right
of the dotted edge). In this procedure,
we have fixed one dotted edge  and thus
multiplied the
contribution of the resulting map
by $-1$
(this is the case 1 of {\bf IP}). The resulting contribution to $\MM$ is
therefore $-\MM(RU_i)\MM(SU_i)$. Thus,
 the first term in \eqref{sanspotentiel} computes the operation of gluing rings by dotted edges.
\end{enumerate}
Putting these  two possibilities together we see that the state $\mu$
 and the enumeration of maps  $\MM$
satisfy the same induction so that they are equal; $\MM(pU_ib)=\mu(pU_ib)$
for any $b$ monomial which does no contain any of the $(U_i, 1\le i\le m$).
Note here
that no dotted edges between rings of incoming arrows can be drawn
since if there are no outgoing arrows in a map, but some $U_i^*$, there is
no contribution.
By traciality of $\mu$, we deduce as well that $\MM_0$ is
tracial. Indeed, if we decompose $p,q$ into $p=p_1U_{i_1}p_2 U_{i_2}\cdots
p_{n-1}U_{i_{n-1}} p_n$ and $q=q_1U_{j_1}q_2 U_{j_2}\cdots
q_{r-1}U_{j_{r-1}} q_r$
with monomials $p_i,q_i$  which  does no contain any of the $(U_i,
1\le i\le m$),
\begin{align*}
\MM(pq)&= \MM( (pq_1U_{j_1}q_2 U_{j_2}\cdots
q_{r-1}U_{j_{r-2}} q_{r-1} U_{j_{r-1}} q_r)\\
&=\mu(pq_1U_{j_1}q_2 U_{j_2}\cdots
q_{r-1}U_{j_{r-2}} q_{r-1} U_{j_{r-1}} q_r)=\mu(pq)\\&
= \mu(qp)=\mu((qp_1U_{i_1}p_2 U_{i_2}\cdots
p_{n-1})U_{i_{n-1}} p_n)=\MM(qp).
\end{align*}

Now we turn to the general $V$ case.

We first check the induction relation when the root star $P$
contains a $U_i$ for some $i\in\{1,\cdots,m\}$ so that we can write
$P=QU_i$. Let us denote for $n$-tuples $\bk=(k_1,\cdots,k_n)$ and
$\bl=(l_1,\cdots,l_n)$, $\binom{\bk}{\bl}=\prod_i\binom{k_i}{l_i}$.
We check the formal equality by considering the induction relation,
now given by:
\begin{align}
\mu^{\bf k+1_j} (QU_i)&=-\sum_{\bl\le \bf k+1_j}
\sum_{Q=RU_iS}\binom{\bk+1_j}{\bl}
\mu^{\bl}(RU_i)\otimes \mu^{\bf k+1_j-\bl}
(SU_i)\nonumber\\&
+\sum_{\bl\le \bf k+1_j}\sum_{Q=RU_i^*S}\binom{\bk+1_j}{\bl}
\mu^{\bl}(R)\otimes \mu^{\bf k+1_j-\bl}(S)\label{munum}\\&
-\sum_{q_j=RU_iS}k_j
\mu^{\bf k}(QU_iSRU_i)-\sum_{q_j=RU_i^*S}k_j
\mu^{\bf k}(QSR)\nonumber\\
\nonumber
\end{align}
We need to show that the enumeration of maps satisfies the same relation. We start by putting stars of type $(q_j, 1\le j\le
n)$ inside a root star
of type  $QU_i$ and we wonder what happens to the root element
$U_i$. We apply
one step of  IP.
Two things can happen. Either we link $U_i$ to another part of $Q$ and in that case we have already shown that the possibilities are enumerated by the
first two terms of the induction relation. Here, note
that the product of $\binom{k_i}{\ell_i}$ corresponds to the possible
distribution of stars
in each part (or submap) of the map, since all the stars
 are labeled.

Thus we need to show that the two other terms take into account the case where $U_i$ is linked to another star of type $q_j$. According to our construction rules we have two possibilities:
\begin{enumerate}
\item
Starting from $U_i$ we glue the arrow to an arrow of the same color entering a star of type $q$. This rule forbids any other gluing from $U_i$, this is counted by
$$\sum_{q_j=RU_i^*S}k_j
\mu(QSR).$$
The coefficient $k_j$ counts the number of choices for the star
 of type $q_j$ since they are all labelled.
\item
The other possibility is to
glue the ring to a ring of the same color. This leads to
$$-\sum_{q_j=RU_iS}k_j
\mu(QU_iSRU_i)$$
possibilities.
\end{enumerate}
In the case where $P$ does not contain any $U_i$, $1\le i\le
m$ but still some $U_i^*$,  the root of the root star
can only be glued by a dotted edge to any other $U_i^*$,
or by a directed edge to a $U_i$ of a star. The resulting
induction relation is exactly given by
 the formula obtained by conjugation of
\eqref{munum}, hence again $M_{\bf k}(P)=\mu^{\bf k} (P)$.
This completes the proof.
\end{dem}

\goodbreak

This theorem gives a combinatorial interpretation in term of maps to
the unitary integrals.
The fact that we do not take the sum
on all maps but only on admissible ones makes this interpretation less
transparent than the
 one for the gaussian case found in \cite{BIPZ}. However, now that we know that the series can be identified to the matrix integral, we obtain  some combinatorial identities
which show that {\bf IP} is less rigid than it looks like.

\begin{corollary}\label{corcombi}
Let $V=\sum t_i q_i$.
\begin{enumerate}
\item
For all $P,Q$,
$$\MM_{\bf t}(PQ)=\MM_{\bf t}(QP).$$
\item
For all monomials $r_1,\dots,r_n,r_{n+1}$, and all permutation $\sigma$ of $n+1$ elements,
$$\MM_{r_1,\cdots,r_n}(r_{n+1})=\MM_{r_{\sigma(1)},\cdots,r_{\sigma(n)}}(r_{\sigma(n+1)}).$$
\item
Assume that we define another procedure to define the root element of the root star (for example we pick the root element to be the second ring avaible if possible, or we pick a ring at random, or any other choice which may change during  {\bf IP} for the root stars that are created during the procedure when new faces are added). This will change the notion of admissible maps and we can define a new weighted sum $\MM'_{r_1,\cdots,r_n}(P)$ and a new series $\MM'_{\bf t}(P)$ where the sum occurs on these new maps. For all $r_1,\dots,r_n,P$,
$$\MM_{r_1,\cdots,r_n}(P)=\MM'_{r_1,\cdots,r_n}(P)$$
$$\MM_{\bf t}(P)=\MM'_{\bf t}(P).$$
\end{enumerate}
\end{corollary}
Note that due to the definition of admissible maps via the procedure
{\bf IP}, those properties are far from being obvious from a purely
combinatorial point of view. Still they will appear as an easy consequence of the identification with the matrix model.

Obviously different roots lead  to a different procedure {\bf IP}, and
thus potentially to different maps.
It is actually possible to see through examples that
this phenomenon actually happens.

However, it follows from the second point of the corollary
that the choice of the root does not affect the weighted sum. The
first and third points show that the choice of the root element
and of the root star  does not affect
the final series.
We were not able to give a more direct combinatorial proof of that result.

To be more specific on the impact
of the choice of the roots on the maps,
let us call clusters the equivalence class
of rings for the equivalence relation
generated by $a\sim b$ if the ring $a$ is glued to the ring $b$ by a dotted edge.
Changing the choices of the roots will lead to different admissible maps since it will allow different positions for the dotted edges. For example, they were three choices for the starting root in
figure \ref{vertex}. For each of these choices, two of the three maps represented
in figures \ref{map1}, \ref{map2} and \ref{map3} would have been reachable by the inductive construction {\bf IP} but not the third one. The one who is not constructible depends on the choice of the first root. It seems that if the maps are different, nevertheless the clusters are the same and in that simple case, knowing this cluster is sufficient to define the faces created by the dotted edges and thus
the weight of the maps.

\begin{dem}

Changing the root element of a star is the same thing than making a circular permutation of the variable of the associated monomial. The theorem shows that weighted sums are equal to the limit of the empirical measure of the matrix model which are tracial. The first and third items are a direct consequence of this identification.

For the second item, observe that permuting the first $n$ monomials doesn't change the sum by its definition. Thus we only need to show that
$$\MM_{r_1,\cdots,r_n}(P)=\MM_{P,r_2,\cdots,r_n}(r_1).$$
Let us define $V=\sum_iu_ir_i+tP$. We will again use the identification with the matrix model but now we will use the formal version.
The coefficient $\MM_{r_1,\cdots,r_n}(P)$ appears as the coefficient
of the limit tracial power state $\muf$
by  Corollary
\ref{convcor} and Theorem \ref{maincombinatorics}.
More precisely,
$$\MM_{r_1,\cdots,r_n}(P)=
\lim_N\left.\frac{\partial^{n}}{\partial u_1\cdots\partial
u_n}\muf(P)\right|_{u_i=0}.$$
We now use the fact that $\muf$ is the limit coefficientwise of the formal model defined in \eqref{etatf}.
Thus,
\begin{align*}
\MM_{r_1,\cdots,r_n}(P)&=\lim_N
\left.\frac{\partial^n}{\prod_i\partial u_i}\frac{E[\mun(P)e^{ N^2\mun ( V)}]}{E[e^{
N^2\mun(V)}]}\right|_{u_i=0 t=0}\\
 &=\lim_N \left.\frac{\partial^{n+1}}{\partial t\prod_i\partial u_i} \frac{1}{N^2}\ln E[e^{
N^2\mun(V)}]\right|_{u_i=0,t=0}.
\end{align*}
We conclude by noticing that this last expression is symmetric in the monomials $r_1$,\dots,$r_n$,$P$.
\end{dem}

\section{Application to free probability}\label{secfreeproba}
In this section we look at applications of the
 combinatorial results of section
\ref{seccombinatorics} to free probability.

Let us assume that the $U_i$'s are chosen independently according to the Haar measure.
If we define $X_i=U_i^*A_iU_i$ then the $X_i$'s are asymptotically
free (according to a theorem of Voiculescu \cite{vo2}) and with fixed
distribution
$\mu$  uniquely defined
 by the distribution of the $A_i$'s. We are interested in using our setup to compute limits of moments
of these variables or in other word to compute the moments of free variables:
$$\mu(X_{i_1}...X_{i_k}).$$
According to our interpretation this can be computed by looking at the maps above the star of type $X_{i_1}...X_{i_k}$ without any other
stars, in other words  we have to focus on computations of
$\MM(q)=\MM_0(q)$ which turns out
to be equal to  $\mu (q)$ where $\mu$ is the free state product
(see Corollary \ref{asfree})

We are interested in using this method to compute some
non-commutative moments of free variables, in relation with
 Speicher's non-crossing cumulants theory, cf \cite{s1}.

\subsection{One star maps}

For these purposes we need to find a simplified interpretation of
$\MM(q)$ in the single star map.

\begin{figure}[ht!]
\begin{center}\resizebox{!}{4cm}{\includegraphics{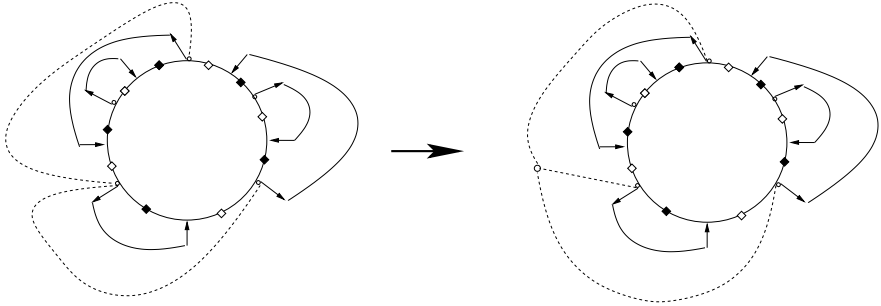}}
\end{center}
\caption{Reduction of a map on one star to a one-star map.}\label{newinterpretation}
\end{figure}

For this case with only one star, the combinatorial interpretation can be slightly modified. First,
we do not need to consider dotted edges between
incoming arrows since if
 there is a $U^*_i$  there must be
 a $U_i$ which can be
 chosen
 as the root element or we can not build any map. But the main difference is that now each time we glue two rings, the edge newly created separate these two rings into two different faces so that they can no longer be glued together. Thus, we can forget about the restriction of the construction rules and present a simpler description in that case. Instead of gluing the ring two by two we will now glue them together. We define a new structure
which we will call a node and now rings can only be glued to node and a node can be glued to any number of rings. A one star map is a map with one star where the arrows has been glued two by two while respecting the orientation and rings may be glued to exactly one node, each node is glued to an arbitrary number of rings but at least one.
Figure \ref{newinterpretation} shows the new representation of a one star map. The trick to go from the previous interpretation to this one is to glue together to a node all the rings that are in the same class of the equivalence relation generated by being glued.
In order to compute
the weight of such a map, observe that several maps give the same one-star map,
 but the weight is easy to compute since as we will see we only need to add a factor $C_{d-1}$ for each node of degree $d$.

\begin{definition}

A one-star map is a connected graph
embedded on a sphere above one star and with some edges
such that
\begin{enumerate}
\item
Edges are drawn only between rings and must not intersect.
\item
Arrows must be glued two by two while respecting the orientation and the color: an arrow going out of a star (associated with a variable $U_i$) is always glued to exactly one other arrow going into a star (associated to a variable $U_i^*$) of the same color. This pair of arrows creates an oriented edge.
\item
Any number of rings may be glued together on a node.
\end{enumerate}
\end{definition}

The weight of a one-star map is the product of the weight of its faces which is defined as before as trace of products of $A_i$'s times the product of the weight of the nodes. The weight of a node of degree $d$ is $(-1)^{d-1}C_{d-1}$. We define
$\TM_0(q)$ the weighted sum of one-star map above a star of type $q$. Note that we no longer need to take care of roots and
of maps that can be built with some set of rules.

\begin{proposition}\label{propappfreeproba}
For all monomial $q$,
$$\mu(q)=\TM_0(q).$$
\end{proposition}

\begin{dem}
We only need to show that $\MM(q)=\TM_0(q)$.
For this we need to compute the number of maps above one star that are reduced
to a given one-star map. The reduction goes as follows: two rings are glued to the same node if they are linked by a sequence of dotted-edges.
We have to count how many configurations of dotted edges lead
to a node of degree $d$. When one of the ring
glued to this node becomes the root in the recursive construction,
it has to be glued to one of the
 other ring glued to the node. Thus it separates the set of ring into two subsets, so according
 to our inductive procedure of
 section \ref{seccombinatorics}, we have to continue to glue this ring to other ones while we continue the construction in the face newly created.
 This yields a structure of tree on this set of rings. We have as many choices as they are trees with $d-1$ edges (to glue the $d$ ring we need exactly $d-1$ edges). This explains the factor $C_{d-1}$. The factor $(-1)^{d-1}$
 simply comes from the factor $-1$ which comes with each edge.
\end{dem}

\subsection{Maps and cumulants}\label{applicfreep}

Let $A_1,\ldots ,A_n$ be self-adjoint variables and $U$ a
unitary matrix,  free from the $A_i$'s.
Then choosing $k$ indices $i_1,\ldots ,i_k$ in $\{1,n\}$ one has
\bigskip
$$\mu(A_{i_1}\ldots A_{i_k})=\mu (U^*UA_{i_1}\ldots U^*UA_{i_k})$$
Let us apply Schwinger-Dyson equation with respect to $U$ to the above equality,
and let us rearrange the sum according to the non-crossing partition
of $A_i$'s generated by the  oriented edges.
Obviously one obtains
a formula of type
\begin{equation}\label{decmom}
\mu(A_{i_1}\ldots A_{i_k})=\sum_{\pi\in NC(k)} \tilde{K}_{\pi}(A_{i_1},\ldots ,A_{i_k})
\end{equation}
where $NC(k)$ is the non-crossing partitions and $\tilde{K}_{\pi}$ is
a $k$-linear form multiplicative along the blocks of $\pi$ in the sense of Speicher: if $\pi=\{V_1,\dots,V_n\}$ with the block $V_i=\{a_1^i,\dots,a_{r_i}^i\}$
$$\tilde{K}_{\pi}(X_{1}\ldots X_{k})=\prod_i{\tilde K}_{(r_i)}(X_{a_{1}^i},\dots,X_{a_{r_i}^i})$$
where $(r_i)$ represents the partition on $r_i$ elements with only one block.

The fact that such a formula holds true for any choice of non-commutative
laws for $A_i$'s proves via the moment-cumulant formula that
$\tilde{K}_{\pi}$ has to be Speicher's non-crossing cumulants $K_\pi$.
But it is also given as a sum on maps by our graphical model.

Let us recap this in the following proposition:

\begin{proposition}
The $n$-th non-crossing cumulant of the variables $A_1$, \dots, $A_p$ is the weight of all one-star maps over
the star build by putting in the clockwise order a ring, a
diamond
of color $i_1$, a ring, a diamond of color $i_2$,\dots, a ring, a diamond of color $i_p$.
\end{proposition}
Note that we have defined this map above a star which is not of type $q$ for any monomial $q$.
This would be a problem for admissible maps since {\bf IP} requires the presence of oriented edges.
But the definition of one-star map is fine in this context.

Actually, Proposition \ref{propappfreeproba} gives us a new proof of the
following Corollary, due to Speicher and known as
non-crossing Moebius formula
\begin{corollary}
The following inversion formula holds true:
$$K_n(A_1,\ldots,A_n)=\sum_{\pi\in NC(n)}\mu_{\pi}(A_1,\ldots ,A_n)
(-1)^{n-|blocks (\pi^c)|}\hspace{-5pt}\prod_{B\,\, block \,\, of \pi^c}
C_{|B|-1},$$
where $\pi^c$ is the Kreweras complement (see \cite{NS})
 and $C_q$ the catalan number.
\end{corollary}
\begin{dem}
This is a direct consequence of the previous proposition. Remember
 that $K_n(A_1,\ldots,A_n)$ is a weighted sum over maps with dotted edges since the star contains some rings and no arrows. These dotted edges form a non-crossing partition of $[|1,\dots,n|]$ by saying that two rings are in the same component if their are linked to a same node. The weight associated to this map is a product whose factors are: $(-1)^{d-1}C_{d-1}$ for each node of degree $d$ and the weight of each face. The faces are by definition the component of the Kreweras complement of $\pi'$.
Thus we obtain:
$$K_n(A_1,\ldots,A_n)=\sum_{\pi'\in NC(n)}\mu_{(\pi')^c}(A_1,\ldots ,A_n)\prod_{B\,\, block \,\, of \pi'}(-1)^{|B|-1}
C_{|B|-1}.$$
The formula follows after taking $\pi'=\pi^c$.
\end{dem}

As a further remark, one can also read graphically the main properties of cumulants, for example,
$K_n(X_{1},\ldots ,X_{n})=0$ as soon as there are occurence of free elements.
More precisely, assume that we can partition the $X_i$'s into two families the $A_j$'s and the $B_k$'s with the algebra generated by
the $A_j$'s free from the algebra generated by the $B_k$'s.
Then if all the $X_i$'s do not take value in the same algebra, $K_n(X_{1},\ldots ,X_{n})=0$.
Indeed, one can replace all the family of $A_j$'s
by the one of $V^*A_jV$ with $V$ unitary and
free from the other variables. Now when looking at the combinatorial interpretation of
$\mu(X_{1},\ldots ,X_{n})$ we can see that the oriented edges coming
from $V$ separate the components
containing the $A_j$'s from the others. By following those edges we see that the faces they are defining contain only variable
from one of the two algebras (The edges are going in the clockwise order around the faces which contain the $B_k$'s and in the counter-clocwise
order around the faces containing the $A_j$'s).
Thus, in the decomposition \eqref{decmom}, the terms corresponding to
partitions with one component containing both some $A_i$'s and some
$B_j$'s vanish.
By uniqueness of the decomposition into cumulants we deduce that those elements vanish i.e. $K_n(X_{1},\ldots ,X_{n})=0$.

These remarks are not new but this shows that our graphical model
fully encompasses the theory of non-crossing cumulants and
that the Schwinger-Dyson equation can also be read in terms of cumulants.

It is interesting to mention here that papers
\cite{HigherOrderFreeness2}
and  \cite{mst}
have developed a calculus on annuli which seems to be
related to our graphical model.
However these approaches only deal with the
asymptotics of second order cumulants whereas
our approach via formal calculus, see section \ref{formalanalyticity},
allows us to deal with arbitrary order
cumulants.

The actual relation can be found in
\cite{cmss}, where
convolution on partitioned permutations is introduced and showed to be
the relevant algebraic tool to handle higher order freeness, namely, the
asymptotic behaviour of cumulants of unitarily invariant random matrices.

 But
the results in our paper give an
explicit algorithmic description of the Moebius inversion formula and therefore
of higher order cumulants. As in the one star case, cumulants are also obtained by
inserting an outer $U^*U$ between each variable of each star
and by looking at generating function where $U$ is linked to its neighboring
$U^*$.

It is interesting to see that a direct (yet difficult to describe) graphical reading of the
Schwinger-Dyson equation, which is our main tool of investigation of unitarily
invariant matrix models, yields non-crossing and could yield
higher order moments related series and operations similar to convolution,
although these latter results rely on more representation theoretic grounds
(Weingarten function theory as developed in \cite{cs}).

It is not obvious to us how  the Schwinger-Dyson
equation can be read off
from the
results of \cite{cmss} (without writing a change of variable invariance formula),
and it would be interesting to attempt to figure out the meaning of Schwinger-Dyson
equation at the representation theoretic level.

\section{Application to the asymptotics of $I_N(V,A_i^N)$}\label{secmatrixintegrals}
In this section, we investigate the free energy by using the combinatorial interpretation of the previous section.

Let $(q_1,\cdots,q_n)$
be fixed monomials in $\cxm$,
let $V=\sum t_i q_i$ be a self-adjoint polynomial
and $I_N(V,A_i)$ be given by \eqref{unitarymatrixint}.

\begin{theorem}\label{mainsec4}
There exists $\varepsilon=\varepsilon(q_1,\cdots,q_n)$
so that for any $\bt\in\C^n\cap B(0,\varepsilon)$
such that $V=V^*$ for any $\alpha\in [-1,1]$,
$$F_{V,\tau}(\alpha):=\hspace{-5pt}\lim_{N\ra \infty}\frac{1}{N^2}\log
I_N(\alpha V_{\bt}, A_i^N) =\hspace{-9pt}\sum_{{\bk}\in \N^n\backslash (0,..,0)}\hspace{-2pt}
\prod_{i=1}^n \frac{(\alpha t_i)^{k_i}}{k_i!}
{\mathbb M}_\bk(q_1,\cdots,q_n,\tau) .$$
Moreover,
$${\mathbb M}_\bk(q_1,\cdots,q_n,\tau)=\sum_{m\mbox{admissible maps
with
} k_i \mbox{ stars } q_i}M_m(\tau)$$
 is the weighted sum of maps constructed above $k_i$ stars of type $q_i$
 for all $i$, after choosing one of them as a root star (this is well defined according to Corollary \ref{corcombi}).
\end{theorem}

\begin{dem}
Let $$F^N_\bt=\frac{1}{N^2}\log
I_N(V_{\bt}, A_i^N).$$
Then, if $\alpha\in\R$,
$$\partial_\alpha F^N_{\alpha \bt}
=\int \mun(V_{\bt}) d\mu^N_{V_{\alpha\bt}}.$$
Assume that $\bt$ is small enough
so that Corollary \ref{convcor} holds
and
remark that $V_{\alpha \bt}$  is self-adjoint
and  such that $|\alpha t_i|\le \varepsilon$
for all $i$ and all $0\le \alpha \le 1$. Thus,
for $\alpha\in [0,1]$,
$$\lim_{N\ra\infty} \partial_\alpha F^N_{\alpha \bt} =
\mu_{\alpha \bt }(V_{\bt})$$
with $\mu_{\alpha \bt }$ the solution
to {\bf SD}$[\alpha V_{\bt}, \tau]$.
By dominated convergence theorem (since $\partial_\alpha F^N_{\alpha
\bt}$
is uniformly bounded in $N$ and $\alpha\in [0,1]$),
we deduce that
$$\lim_{N\ra\infty}  F^N_{\alpha \bt} =\int_0^1
\mu_{\alpha \bt }(V_{\bt})d\alpha $$
where we used that $F^N_0=0$.

\end{dem}

Here also, we obtain the following important corollary, as a consequence
of Corollary \ref{corSDf}.

\begin{corollary}
The following holds true:
\begin{align*}\lim_{N\ra\infty}
& \frac{\partial ^k}{\partial z^k}N^{-2}\log \int_{\UNCm} e^{zN Tr
 (V(U_i,U_i^*,A_i^N, 1\le i\le m ))} dU_1\cdots dU_m|_{z=0}\\
 &=
\frac{\partial ^k}{\partial z^k} F_{V,\tau}(z)|_{z=0}\end{align*}
\end{corollary}

In particular, this result allows us to give an expansion of the
Harish-Chandra-Itzykson-Zuber integral as a
generating function
of the number of some  maps.
Let us recall the exact expression of this integral:
$$F_N^{A,B}(z):=\frac{1}{N^2}\log HCIZ(zA,B)=
\frac{1}{N^2}\log\int_{\Ua_N}e^{zNTr(U^*AUB)}dU.$$
The maps appearing in the expansion contain only stars
of type  $U^*AUB$ (see the star in the middle of figure \ref{vertex}).
Besides we can build these maps without considering the rings attached
to variable $U^*$ since we will always be able to choose the root
element to
be  a $U$ (a $U^*$ always comes with a $U$ for this potential).

Since the number of diagrams
is growing quickly we compute
only the first term of the expansion.
Note that when gluing the arrow of the  root of the root star,
we must always glue it to another incoming arrow
of another  star and hence we shall never see the case
of a root star with no $U_i$'s. Again, we therefore do not see
dotted edges between incoming arrows.

 Besides, we consider only the case where the distribution is centered,
 that is when $\tau(A)=\tau(B)=0$.
The other cases can be deduced easily from this one since we have the relation
$$F_N^{a+A,b+B}(z)=F_N^{A,B}(z)+\frac{z}{N}(b\Tr A+a\Tr B) + zab.$$
In terms of diagrams, this means that we only need to consider diagrams
 such that no face contains only one diamond.

According to the previous theorem, $\lim_{N\ra\infty}F_N^{A,B}(z)$ has, for small $z$, an expansion $\sum_n F_n z^n$.
We now use this graphical representation to compute the first terms of this integral.

Since the distributions are centered, the first
term
$F_1$ is zero.

\begin{figure}[ht!]
\begin{center}\resizebox{!}{4cm}{\includegraphics{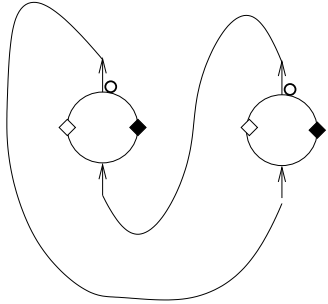}}
\end{center}
\caption{Second term in the expansion of the HCIZ integral.}\label{iz2}
\end{figure}

The second term $F_2$ consists
of maps constructed with two stars of type $U^*AUB$. There is only one way to add edges between these two stars
to construct a connected map without faces which contains only one diamond, this is represented by figure \ref{iz2}.
We obtain a map with two faces.
One has two diamonds associated to $A$ and the other one two
diamonds
associated to $B$. Thus the weight of this map is $\tau(A^2)\tau(B^2)$. Since there is no gluing between the rings they are no other signs. They are only one way to distribute the labels on this picture (that is the second distribution leads to the same map) thus to obtain $F_2$ we only need to divide by $2!$,
$$F_2=\frac{1}{2}\tau(A^2)\tau(B^2).$$

\begin{figure}[ht!]
\begin{center}\resizebox{!}{4cm}{\includegraphics{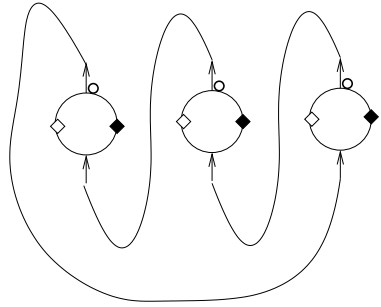}}
\end{center}
\caption{Third term in the expansion of the HCIZ integral.}\label{iz3}
\end{figure}

We can continue this for the next terms in the expansion, the third term (see figure \ref{iz3}) is in the same spirit and leads to
$$F^3=\frac{1}{3}\tau(A^3)\tau(B^3).$$

\begin{figure}[ht!]
\begin{center}\resizebox{!}{10cm}{\includegraphics{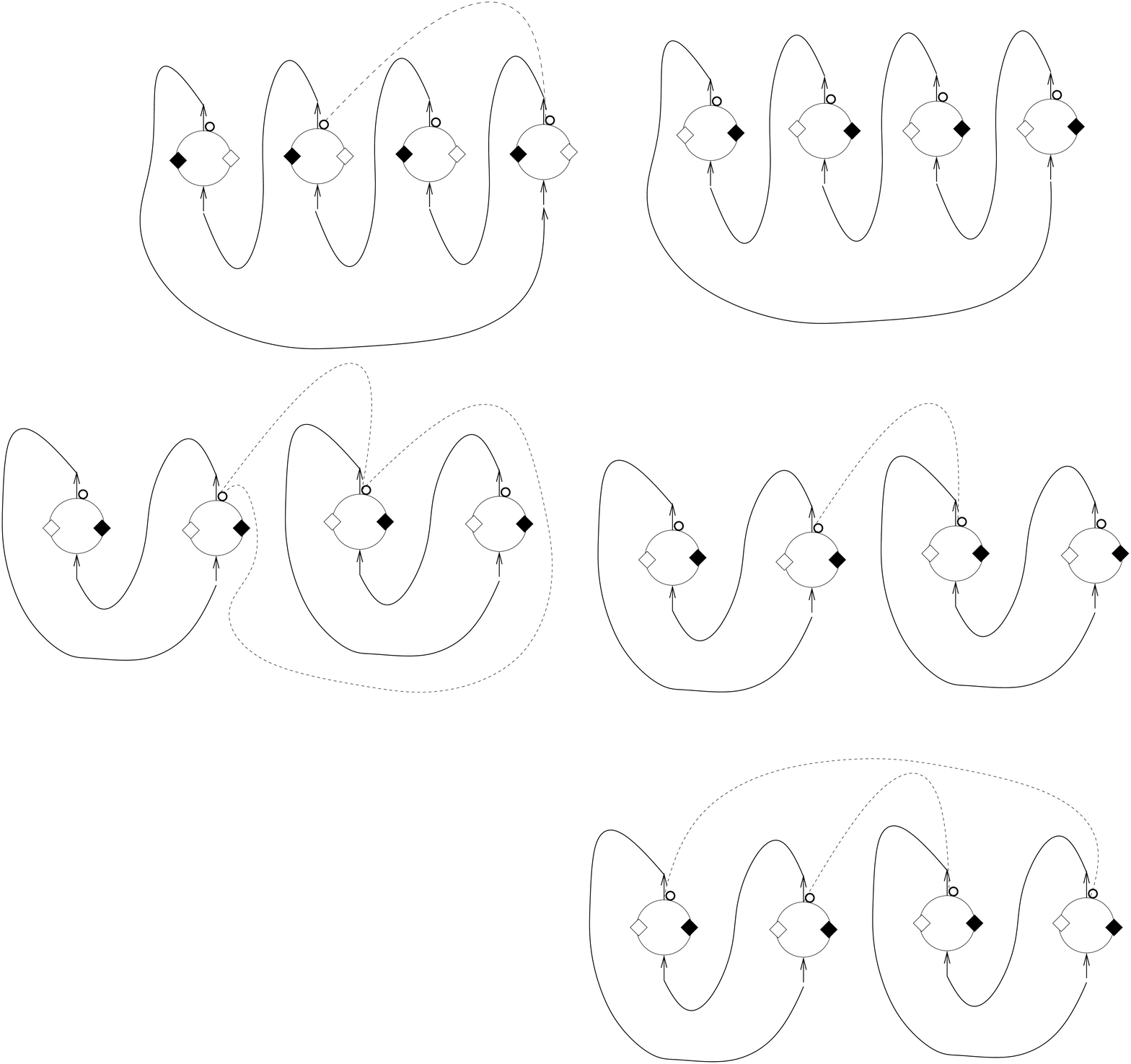}}
\end{center}
\caption{Fourth term in the expansion of the IZ integral.}\label{iz4}
\end{figure}

The fourth term is the first one where gluings between the
 rings appear. Thus weigths with negative coefficients can occur. The sign of a map is easy to compute, it is $-1$ to the power the number of dotted lines in the map. Equivalently since in the
case of HCIZ integral the number of oriented edges
 is equal to the number of stars, this number is also equal to the number of faces of the map and thus to the number of factor in the product of moments of the weight. In figure \ref{iz4}, we have drawn all unlabelled planar maps one can construct with
$4$ stars. To compute the exact coefficient of each map one has to multiply it by the number of way to distribute the
labels and divide by $4!$.

This leads to,
$$F^4=\frac{1}{4}\tau(A^4)\tau(B^4)-\frac{1}{2}\tau(A^2)^2\tau(B^4)-\frac{1}{2}\tau(A^4)\tau(B^2)^2$$
$$+\frac{1}{2}\tau(A^2)^2\tau(B^2)^2+\frac{1}{4}\tau(A^2)^2\tau(B^2)^2.$$
Here the weight are given in the same order than the maps in the figure.
Note a new and interesting feature that appears in the third map: two rings are linked by more than one
dotted edge.

\begin{figure}[ht!]
\begin{center}\resizebox{!}{10cm}{\includegraphics{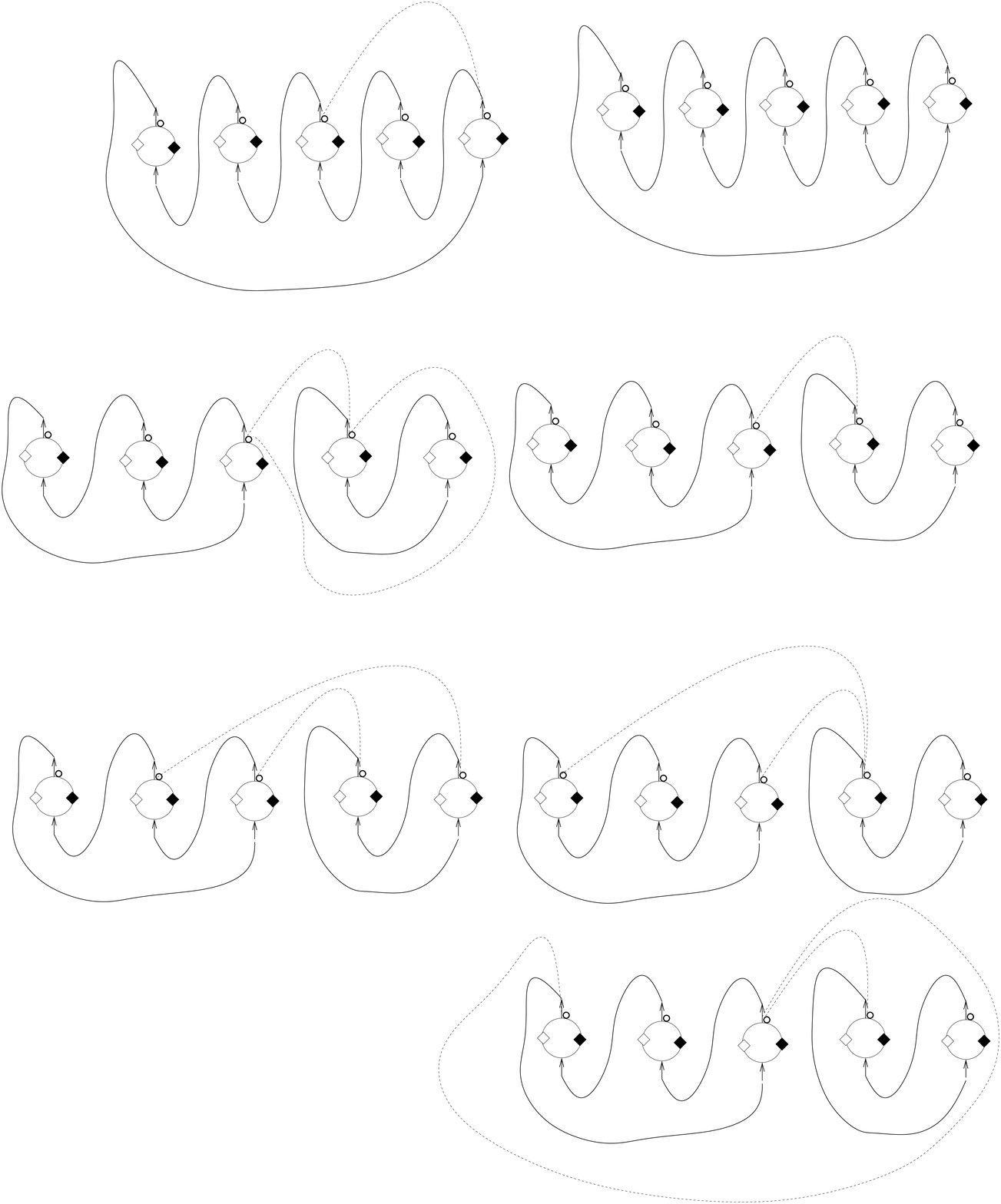}}
\end{center}
\caption{Fifth term in the expansion of the IZ integral.}\label{iz5}
\end{figure}

The other terms can be computed in the same way, for example figure \ref{iz5} represents the fifth term and gives
$$F_5=\frac{1}{5}\tau(A^5)\tau(B^5)-\tau(A^2)\tau(A^3)\tau(B^5)-\tau(A^5)\tau(B^2)\tau(B^3)$$ $$+4\tau(A^2)\tau(A^3)\tau(B^2)\tau(B^3).$$

Thus the first terms agree with the expansion given in
\cite{ZZ} on page 23, besides this allows us to answer a question
raised in this paper. Indeed, the
authors ask if there is an explanation to the fact that the coefficient of
$F_n$ all seem to be integer multiple of $\frac{1}{n}$. This is easy to prove with this graphical interpretation. To compute
the contribution of a given unlabelled map we must distribute the labels $\{1,..,n\}$ on its stars, count the number of different map that we obtain and divide by $n!$. But after choosing the star which received the label $1$ we have $(n-1)!$
ways to distribute the remaining labels and they all lead to different maps (note that on the other hand, due
to possible symmetry in the unlabelled map, different choices for the star with the label $1$ may lead to the same maps). Thus the coefficient in front of this map is a multiple of $\frac{(n-1)!}{n!}=1/n$. More precisely it is $1/n$ times the number of choices of the star which carry the label $1$ that will lead to different maps, in particular it is always less than $1$.

To finish, we wish to point out that we can recover results in \cite{c1}
and \cite{GM1} about scalings of HCIZ integral.
In these two papers, one considers the
scaling where $A$ has small rank, which amounts
to considering only
terms $\tau(A^k)\times P(B)$. Here the transformation depicted in section
\ref{secfreeproba}
applies and one sees that $P(B)$ has to be $k^{-1}K_k(B)$.
In particular this means in the case that $A$ is a rank $1$ projection, that
$N^{-1}\log HCIZ$ tends to the primitive of Voiculescu's $R$-transform.

\section{Application to Voiculescu free entropy}\label{secfreeentropy}
Voiculescu's microstates free entropy
is given as the asymptotic
 the volume of matrices whose
empirical distribution
approximates sufficiently well
a given tracial state. Up to a Gaussian
factor, it is given
by
$$\chi(\mu)=\limsup_{\sur{\e\downarrow 0}{
k\uparrow \infty ,R\uparrow \infty}}\limsup_{N\ra\infty}
\frac{1}{N^2}\log \mu_N^{\ot m}\left(\Gamma_R(\mu,\e,k)\right)$$
with $\mu_N$ the Gaussian measure on $\Ha_N$ and $\Gamma_R(\mu,\e,k)$
the microstates
$$\Gamma_R(\mu,\e,k)=\{X_1,\cdots,X_m\in\Ha_N:
|\frac{1}{N}\Tr(X_{i_1}\cdots X_{i_p})-\mu(X_{i_1}\cdots X_{i_p})|<\e$$
$$\qquad\qquad p\le k, i_\ell\in\{1,\cdots, m\}, \|X_i\|_\infty\le R
\}.$$
When $m=1$, it is well known \cite{vo}
that $\mu\in\Pa(\R)$ and

$$\chi(\mu)=I(\mu)=\int\int\log|x-y|d\mu(x)d\mu(y)
-\frac{1}{2}\int x^2 d\mu(x) +\mbox{const.}$$
Moreover, one can replace
the $\limsup$ by a $\liminf$ in the definition
of $\chi$. Such answers (convergence and
formula for $\chi$) are still open in general  when $m\ge 2$
(see \cite{BCG} for bounds). However, if $\mu$ is the  law of
$m$ free variables with respective laws $\mu_i$,
then these questions are settled
and
$$\chi(\mu)=\sum_{i=1}^m I(\mu_i).$$
 We here want to
emphasize that our result
provides a small step towards
dependent variables by showing convergence and
giving a formula for the type
of laws $\mu$ solutions of Schwinger-Dyson's equations ${\bf SD[V,\tau]}$.
Indeed, we shall prove that
\begin{theorem}\label{thmfreeentropy}
Let $\mu$ be the law
of $m$ self-adjoint variables $X_i$
with marginal distribution $(\mu_1,\cdots,\mu_m)$.
Assume that $X_i$ can be decomposed as
$X_i=U_iD_iU_i^*$ with $U_i$ unitary matrices
in such a way that the joint law $\nu$ of $(D_i,U_i,U_i^*)_{1\le i\le m}$
satisfy  ${\bf SD[V,\tau]}$ with $\tau$ the law
of $m$ free variables with marginal
distribution $\mu_1,\cdots,\mu_m$ and some potential $V=\sum_{i=1}^n t_i q_i$.
Assume  that the $t_i$'s are small enough
so that Corollary \ref{convcor} holds.
Assume also that the hypotheses of
Theorem \ref{mainsec4} hold.
Then,
$$\chi(\mu)=\liminf_{\sur{\e\downarrow 0}{
k\uparrow \infty}}\liminf_{N\ra\infty}
\frac{1}{N^2}\log \mu_N^{\ot m}\left(\Gamma_R(\mu,\e,k)\right)$$
and a formula of $\chi(\mu)$ can be given in terms of the $\mu^{\bf
k}$'s
of Theorem \ref{regularity}.
\end{theorem}

\begin{dem}
Indeed, let us  consider
$V=V(U_iA_i U_i^*, 1\le i\le m)$
with $V$ a self-adjoint
polynomial
and $\mu$ the unique solution
of ${\bf SD[V,\tau]}$
with $\tau$ the law of the $A_i,1\le i\le m$
which is now chosen to be the law of $m$ free
variables with marginals distribution $\mu_i$, $1\le
i\le m$.
Under the law $\mu_N^{\ot m}$, we can diagonalize the
matrices $X_i=U_i D_iU_i^*$
with $U_i$ following the Haar measure on
$\UNC$, and id $d$ is the Dudley metric,
we find that for $N$ sufficiently large
\begin{eqnarray*}
\L_N&:=&\mu_N^{\otimes m}
\left( \Gamma_R(\mu,\e,k)\right)\\
&=&\mu_N^{\otimes m}
\left(
d(\mun_{D_i}, \mu_{i})<\e;\mun_{ U_iD_iU_i^*,1\le i\le m}
\in\Gamma_R(\mu,\e,k)\right)
\\
&=&\int_{\sur{d(\mun_{D_i}, \mu_{i})<\e}{\|D_i\|_\infty\le R}}
  \left( \int_{\mun_{ (U_iD_iU_i^*)_{1\le i\le m}}
\in\Gamma_R(\mu,\e,k)} dU_1\cdots dU_m\right)\prod_{1\le i\le m}
d\sigma_N({\bf\l}_i
)
\\
\end{eqnarray*}
where we denoted $\D(\l_j)=\prod_{k\neq j} |\lambda_k-\lambda_j|$
and $d\sigma_N$ the probability measure
$$d\sigma_N({\bf\l}):=Z_N^{-1} \prod_{k\neq j} |\lambda_k-\lambda_j|^2
 e^{-\frac{N}{2}
\sum (\l^j)^2} \prod_{1\le j\le N} d\l^j.$$
In these notations, $D_i=\mbox{diag}(\l_1^i,\cdots,\l_N^i)$
and ${\bf\l}=(\l_1,\cdots,\l_N)$.
Hereafter, $\hat\mu^N_{\{E_i\}_{1\le i\le n}}$
denotes the empirical ditribution of $\{E_i\}_{1\le i\le n}$;
$\hat\mu^N_{\{E_i\}_{1\le i\le n}}(P)=N^{-1}\Tr(P(E_i,{1\le i\le
n}))$.
As a consequence, applying the large deviations result of \cite{BAG}
to the diagonal matrices $D_i$, we find that there exists $o(1)$
 going
to zero with $\e$ such that
\begin{align*}
\L_N&\le e^{N^2\sum_{i=1}^m I(\mu_{i})+N^2 o(1)}
\sup_{\sur{d(\mun_{D_i}, \mu_{i})<\e}{\|D_i\|_\infty\le R}} \int_{\mun_{ \{U_iD_iU_i^*\}_{1\le i\le m} }
\in\Gamma_R(\mu,\e,k)}  dU_1\cdots dU_m\\
&:= e^{N^2\sum_{i=1}^m I(\mu_{i})+N^2 o(1)}  \L_N^1\\
\end{align*}
with for $k$ greater than the degree of $V$,
\begin{align*}
\L_N^1&=
\sup_{\sur{d(\mun_{D_i}, \mu_{i})<\e}{ \|D_i\|_\infty\le R}}
\int_{ \mun_{ \{U_iD_iU_i^*\}_{1\le i\le m} }\in\Gamma_R(\mu,\e,k)}
e^{N\Tr(V)-N\Tr(V)}
dU_1\cdots dU_m\\
&= e^{-N^2\mu(V)+N^2 \e}
\sup_{\sur{d(\mun_{D_i}, \mu_{i})<\e}{ \|D_i\|_\infty\le R}}\int_{
\mun_{ \{U_iD_iU_i^*\}_{1\le i\le m} } \in\Gamma_R(\mu,\e,k)}
e^{N\Tr(V)}
dU_1\cdots dU_m\\
&\le e^{-N^2\mu(V)+N^2 \e}
\sup_{\sur{d(\mun_{D_i}, \mu_{i})<\e}{ \|D_i\|_\infty\le R}}\int
e^{N\Tr(V)}
dU_1\cdots dU_m\\
&=  e^{-N^2\mu(V)+N^2 \e}
\sup_{d(\mun_{D_i}, \mu_{i})<\e,\|D_i\|_\infty\le R} I_N(V,D_i)\\
\end{align*}
Now, for fixed $R$, any $D_i,D_i'$ in
$d(\mun_{D_i}, \mu_{i})<\e,\|D_i\|_\infty\le R$
$$\left|\frac{1}{N^2}\log I_N(V,D_i)-\frac{1}{N^2}\log I_N(V,D_i')
\right|\le\eta(\e,R), $$
with $\eta(\e,R)$ going to zero as $\e$ goes to zero for any fixed
$R$. Hence,
$$\limsup_{N\ra\infty} \frac{1}{N^2}\log I_N(V,D_i)\le
F(V,\mu_i)+\eta(\e,R)$$
with $F(V,\mu_i)$  the limit of $N^{-2}\log I_N(V,A_i)$
given in Theorem \ref{mainsec4} when the distribution of the $A_i$
converges to free variables with marginal
distribution $\mu_i$.
We thus have proved, letting $\e$ going to zero and then $R,k$
to infinity, that

$$\chi(\mu)\le \sum_{i=1}^m I(\mu_{A_i})-\mu(V)
+F(V,\mu_i).$$
Conversely,  we have
\begin{eqnarray*}
\L_N&\ge &e^{N^2\sum_{i=1}^m I(\mu_{i})+N^2 o(\e)} \L_N^2\\
\end{eqnarray*}
with
\begin{eqnarray*}
\L_N^2&:= &
\inf_{d(\mun_{D_i}, \mu_{i})<\e,\|D_i\|_\infty\le R} \int_{\mun_{ (U_iD_iU_i^*)_{1\le i\le m}}
\in\Gamma_R(\mu,\e,k)}  dU_1\cdots dU_m\\
&=& e^{-N^2\mu(V)+N^2 o(\e)}
\inf_{\sur{d(\mun_{D_i}, \mu_{i})<\e}{ \|D_i\|_\infty\le R}}
\int_{ \mun_{ \{
U_iD_iU_i^*\}_{1\le i\le m}}\in\Gamma_R(\mu,\e,k)}
e^{N\Tr(V)}
dU_1\cdots dU_m\\
&\ge & e^{-N^2\mu(V)+N^2 o(\e)}
\inf_{\sur{d(\mun_{D_i}, \mu_{i})<\delta}{ \|D_i\|_\infty\le R}}
\int_{ \mun_{ \{U_iD_iU_i^*\}_{1\le i\le m}}\in\Gamma_R(\mu,\e,k)}
e^{N\Tr(V)}
dU_1\cdots dU_m\\
\end{eqnarray*}
for any $\delta<\e$. Now, choosing $\delta$
and using the continuity of $\mun_{ \{U_iD_iU_i^*\}_{1\le i\le m}}$
in the distribution of the uniformly bounded variables $D_i$,
we find by Corollary \ref{convcor} and
our hypothesis  that
$$\liminf_{N\ra\infty}\frac{ \int_{ \mun_{ U_iD_iU_i^*,1\le i\le m}\in\Gamma_R(\mu,\e,k)}
e^{N\Tr(V)}
dU_1\cdots dU_m}{\int e^{N\Tr(V)}
dU_1\cdots dU_m}=1
$$
which insures that
$$\chi(\mu)\ge\sum_{i=1}^m I(\mu_{i})-\mu(V)
+F(V,\mu_i).$$

Thus we have proved that
$$\chi(\mu)=\sum_{i=1}^m I(\mu_{i})-\mu(V)
+F(V,\mu_i).$$
Note that $\mu(V)$ and $F(V,\mu_i)$
can be written in terms of the $\mu^{\bf
k}$
of Theorem \ref{regularity} by Theorem \ref{mainsec4}.

\end{dem}

\section{Generalization to integrals over the orthogonal group}
In a recent article \cite{Z}, Zuber shows that the large $N$
asymptotics of two matrix integrals (the integral with external
magnetic field and the  Harish-Chandra-Itzykson-Zuber integral)
enjoy a universality property in the sense that they are the same
(up to a proper rescaling) if one integrates over the unitary or the
orthogonal group. This property was also obtained (but not
explicitly stated) in the case of the Harish-Chandra-Itzykson-Zuber
integral in \cite{GZ} where the rate functions for the large
deviation principle for the law of the spectral measure process of
the Hermitian and the symmetric Brownian motion were shown to differ
only by a factor two. The Harish-Chandra-Itzykson-Zuber integral
is  rather special in the family of angular integrals and one can
compute many interesting related quantities, regardless of the group on which
integration is taken (see \cite{BE,FPEDZ}).

In this section, we generalize this universality property by
relating the large $N$ limit of 
any small parameter integrals over the orthogonal group
with its complex analogue.

Let us define
\begin{equation}
I_N^1(V, A_i^N):=\int_{\ONCm} e^{N \Tr (V(O_i,O_i^*,A_i^N, 1\le i\le m))}
 dO_1\cdots dO_m\end{equation}
where $(A_i^N, 1\le i\le m)$ are $N\ts N$ deterministic  symmetric
uniformly bounded  matrices, $dO$ denotes the Haar measure on the
orthogonal group $\ONC$ (normalized so that $\int_{\ONC}dO=1$). In
this section we will assume that $V$ is a non-commutative polynomial
in the $O_i,O_i^*,A_i^N$ with {\it real} coefficients. Here,
$O^*=O^t$ is the standard involution $O^*_{ij}=O_{ji}$. Observe that
if $P$ is a polynomial, $P(O_i,O_i^*,A_i^N, 1\le i\le
m)^t=P^*(O_i,O_i^*,A_i^N, 1\le i\le m)$ so that we keep also the
notation $P^*$.

We then claim that we have the following analogue of Theorem \ref{mainsec4}, which shows that
the first order of integrals over the orthogonal group is the same as on the unitary group (up to proper renormalizations);

\begin{theorem}\label{mainsec4o}
There exists $\varepsilon=\varepsilon(q_1,\cdots,q_n)$
so that for any $\bt\in\R^n\cap B(0,\varepsilon)$
such that $V=V^*=\sum t_i q_i$, if we define
$$F^1_{V,\tau}:=\lim_{N\ra \infty}\frac{1}{N^2}\log I_N^1(V_{\bt}, A_i^N)$$
then $F^1_{V,\tau}$ exists and
$$F^1_{\frac{1}{2}V,\tau}=\frac{1}{2}\sum_{{\bk}\in \N^n\backslash (0,..,0)}
\prod_{1\le i\le n} \frac{t_i^{k_i}}{k_i!}
{\mathbb M}_\bk(q_1,\cdots,q_n,\tau) .$$
Moreover,
$${\mathbb M}_\bk(q_1,\cdots,q_n,\tau)=\sum_{m\mbox{admissible maps
with
} k_i \mbox{ stars } q_i}M_m(\tau)$$
 is the weighted sum of maps constructed above $k_i$ stars of type $q_i$
 for all $i$, after choosing one of them as a root star .
 \end{theorem}
The proof is based on the fact that if   $\mu^{N,1}_V$ denotes
the law on $\ONCm$ given by
$$\mu^{N,1}_{\frac{1}{2}V}(dO_1,\cdots,dO_m):=\frac{1}{I_N^1(\frac{1}{2}V, A_i^N)}
 e^{\frac{N}{2} \Tr (V(O_i,O_i^*,A_i^N, 1\le i\le m))}
 dO_1\cdots dO_m$$
and $\mun$ is the empirical distribution of $(O_i,O_i^*,A_i,1\le i\le m)$,
then we have the analogue of Corollary \ref{convcor}.

\begin{theorem}\label{convcor2}
Assume that $V=\sum t_i q_i$ is self-adjoint. Let $D$ an integer and $\tau$ a tracial state in $\Ma|_{(A_i)_{1\le i\le m}}$  be given.
There exists $\varepsilon=\varepsilon(D,m)>0$
such that if $|t_i|\le \varepsilon$,
$\mun$ converges almost surely under $\mu^{N,1}_{\frac{1}{2}V}$ to the unique
solution $\mu_{\bt}$ of the Schwinger-Dyson equation {\bf SD[V,$\tau$]}.
Moreover, $\bar\mu^{N,1}_{\frac{1}{2}V}=\mu^{N,1}_{\frac{1}{2}V}(\mun)$ converges as
well to this solution
as $N$ goes to infinity.
\end{theorem}
In fact, since then we know that $\mu_{\bt}(P)$ expands
as a generating function of the ${\mathbb M}_\bk(q_1,\cdots,q_n,\tau)$'s,
Theorem \ref{mainsec4o} follows readily  since for any $\alpha\in [0,1]$,
$$\partial_\alpha \frac{1}{N^2}\log I_N^1(\frac{\alpha}{2}V_{\bt}, A_i^N)
=\frac{1}{2} \bar\mu^{N,1}_{\frac{1}{2}V}(V)$$
converges towards $\frac{1}{2} \mu_{\bt}(V)$.

\nn
{\bf Proof of Theorem \ref{convcor2}}
The proof follows
the same lines as the proof of Theorem \ref{mainsec1};
 we make the
change of variables $\bO=(O_1,\cdots,O_m)\in \ONCm
\to \Psi(\bO)=(\Psi_1(\bO),\cdots,\Psi_m(\bO))\in \ONCm$ with
 $$\Psi_j(\bO)= O_je^{\frac{\l}{N} P_j(\bO)}$$
where the $P_j$ are antisymmetric polynomials (i.e. $P_j^*=-P_j$).
The only change is that now $ P_j(\bO)$ are matrices with real
coefficients and the differentials hold in the direction of $\ANR$
which are the antisymmetric matrices with real coefficients. For $N$
large enough, $\Psi$ is a diffeomorphism; it is as in the complex
case a local diffeomorphism which is injective. As such, its image
is open and compact. $ \ONCm$ is not connected but the union of
copies of $SO^\e(N)=\{O\in \ONC;\det(O)=+\e\}$, $\e=+1$ or $-1$.
Since $\det(\Psi_j(\bO))=\det(O_j)\det(e^{\frac{\l}{N}
P_j(\bO)})=\det(O_j)$, $\Psi$ maps  $SO^{\e_1}(N)\ts SO^{\e_2}(N)
\ts\cdots SO^{\e_m}(N)$ into itself for each choice of
$\e_i\in\{1,-1\}$. Therefore, by connectedness of this set,
$\Psi(SO^{\e_1}(N)\ts\cdots \ts SO^{\e_m}(N))$ is open
and closed and therefore equals $SO^{\e_1}(N)\ts SO^{\e_2}(N)
\ts\cdots SO^{\e_m}(N)$. Thus, $\Psi$ is a diffeomorphism of
$\ONCm$. Like in the proof of Lemma \ref{lempsi}, we need to compute
the Jacobian of this change of variable. The same arguments apply to
show that
$$|\det J_\Psi(\bO)|=\exp (\frac{\l}{N}\Tr {\tilde \Phi} + O(1))$$
with ${\tilde \Phi}$ the linear operator defined on antisymmetric matrices by
$${\tilde \Phi}.A = \sum_i\partial_iP_i\sharp A.$$

A basis of $\ANR$ is given, for $k< l$,  by
$$E^1(kl)_{rj}=\frac{1_{r=k,j=l}-1_{r=l,j=k}}{\sqrt 2}.$$
Therefore, the trace of  any linear endomorphism
$\varphi$ on $\ANR$ defined  by $\varphi(X)=\sum_\ell A_\ell XB_\ell$, for uniformly bounded matrices $A_\ell$, $B_\ell$,
is now given by
\begin{eqnarray*}
\Tr(\varphi)&=&\sum_{k<l}\Tr( E^1(kl)^* \varphi(E^1(kl)))=\frac{1}{2}
\sum_\ell(\sum_{k\neq l} A^\ell_{ll}B^{\ell}_{kk} -\sum_{k\neq l}
A_{lk}^\ell B_{lk}^\ell)\\
&=&\frac{1}{2}\sum_\ell \Tr(A^\ell)\Tr(B^\ell) +\Tr (A_\ell B_\ell^t)\\
&=&\frac{1}{2}\sum_\ell  \Tr(A^\ell)\Tr(B^\ell) +NO(1)
\end{eqnarray*}
since the operator norm of $A_\ell$ and $B_\ell$ is uniformly bounded,
$O(1)$ is uniformly bounded
 in $N$.

We can apply this bound to our case where
$A_\ell$ and $B_\ell$ are given by  $\partial_i P_i=:\sum_\ell A_\ell\otimes
B_\ell$. The $A_\ell$ and $B_\ell$'s are uniformly bounded since the  $O_j$'s and the $A_j$'s are
and non zero for a finite number of $\ell$'s, thus
we deduce that
$$|\det J_\Psi(\bO)|=\exp(\frac{\l}{2 N}\sum_{i=1}^m \Tr\otimes\Tr(\partial_i P_i) +O(1))$$
with $O(1)$ bounded uniformly in $N$.
Since $O(1)$ is uniformly bounded, we can now proceed exactly as
in the proof of Theorem \ref{mainsec1} to show that for any
$r\in\{1,\cdots,m\}$,
$$\lim_{N\ra\infty} \left\lbrace\frac{1}{2}\mun\otimes\mun(\partial_r P) +
\frac{1}{2N}\mun(D_rV P)\right\rbrace
=0\hspace{1cm} \mu^{N,1}_{\frac{1}{2}V} \, a.s.$$
As a consequence, for any limit point $\tau$ of $\mun$, any
antisymmetric polynomial $P$,
\begin{equation}\label{SDorth}\tau\otimes\tau(\partial_r P)+\tau(D_rV
P)=0.
\end{equation}
If $P$ is symmetric, we claim that for any $r\in\{1,\cdots,m\}$,
\begin{equation}\label{SDorth2}
\tau\otimes\tau(\partial_r P)=\tau(D_rV P)=0
\end{equation}
so that \eqref{SDorth} still holds. Indeed,  if $Q$ is a word in the
$( O_i,O_i^* ,A_i,1\le i\le m)$,
\begin{align*}
\partial_r Q
&=\sum_{Q=Q_1O_r Q_2} Q_1O_r\otimes Q_2-\sum_{Q=Q_1O^*_r Q_2}
Q_1\otimes O_r^*Q_2\\
\partial_r Q^*&=\sum_{Q^*=Q_1O_r Q_2} Q_1O_r\otimes
Q_2-\sum_{Q^*=Q_1O^*_r Q_2} Q_1\otimes O_r^*Q_2\\
&=\sum_{Q=Q^*_2O^*_r Q^*_1} Q_1O_r\otimes Q_2-\sum_{Q=Q^*_2O_r
Q^*_1} Q_1\otimes O_r^*Q_2\\
&=\sum_{Q=Q_1O^*_r Q_2} (O^*_rQ_2)^*\otimes Q^*_1-\sum_{Q=Q_1O_r
Q_2} Q_2^*\otimes (Q_1O_r)^*.
\end{align*}
Since the trace is invariant under transposition, we deduce that for
all $P$, $\mun(P^*)=\mun(P)$ and thus,
$$\mun\otimes\mun(\partial_r
Q+\partial_r Q^*)=0.$$

With the same method, we can deal with the cyclic derivative term.
Indeed, since $D_r(Q^*)=-(D_r Q)^*$, if we write $V=Q+Q^*$, we
obtain:
\begin{align*}\mun(D_rV
(P+P^*))&=\mun(D_r(Q+Q^*)(P+P^*))\\
&=\mun(D_rQ (P+P^*))-\mun((D_rQ)^*(P+P^*))\\
&=\mun(D_rQ (P+P^*))-\mun((P+P^*)D_rQ)=0.
\end{align*}
To sum up,
$$\mun\otimes\mun(\partial_r P)=\mun(D_rV P)=0$$
from which we get \eqref{SDorth2} by going to the limit. Since any
polynomial $P$ can be decomposed as the sum of a symmetric polynomial
($P+P^*/2$) and an antisymmetric polynomial ($P-P^*/2$), we conclude
by linearity that \eqref{SDorth} holds for any polynomial $P$. By
uniqueness of the solutions to this equation for sufficiently small
parameters $t_i$ proved in Theorem \ref{mainsec2}, the proof is
complete.
\hfill\qed

{\bf Acknowledgments:} During the Spring 2007 (during which the
chapter 5 was worked out) Alice Guionnet
 visited the UC Berkeley
Department of Mathematics. Her visit was supported
in part by funds from NSF Grants  DMS-0405778,
DMS-0605166 and DMS-0079945.
Beno\^\i t Collins' research was partly supported by an NSERC grant. \'Edouard Maurel-Segala visited Stanford University during 2006-007 and was
supported by funds from NSF grant DMS-0244323.We thank J-B. Zuber
for very fruitful discussions, which in particular led us to insert the last section of this article.


\begin{thebibliography}{99}
\bibitem{APS}{Albeverio S.,  Pastur L., Shcherbina M.},
     {On the {$1/n$} expansion for some unitary invariant ensembles
              of random matrices},
       {\it  Comm. Math. Phys.},{\bf{224}},{271--305} (2001)

\bibitem{ACKM}{Ambjorn J.  and Chekhov  and   Kristjansen  and C.F. and
Makeenko Yu.},{Matrix model calculations beyond the spherical limit},
 {\it  Nuclear Physics B},
{\bf{404}},{127-172},({1993})
\bibitem{BAG}Ben Arous G. and  Guionnet A.,
Large deviations for {W}igner's law and {V}oiculescu's
              non-commutative entropy
{\it  Probab. Theory Related Fields},{\bf{108}},{517--542}
({1997}).
\bibitem{BE} Berg\`ere M. and Eynard B., Some properties of angular
integrals, {\em arxiv: 0805.4482}(2008)


\bibitem{BCG}
Biane P.,  Capitaine M. and  Guionnet A.,
Large deviations bounds for the
 law  of the trajectories of the Hermitian
Brownian motion {\it Inventiones
 Mathematicae},
{\bf 152}, 433-459(2003).


\bibitem{Bi}
    Biane P.,
     Logarithmic {S}obolev inequalities, matrix models and free
              entropy,
{\it Acta Math. Sin. (Engl. Ser.)},
{\bf 19}-3, 497--506, (2003).


\bibitem{BIPZ}  Br{\'e}zin E., Itzykson C. , Parisi G. and  Zuber J.
              B., Planar diagrams {\em  Comm. Math. Phys.}{\bf 59},
35--51,(1978).


\bibitem{c1} Collins B., Moments and cumulants of polynomial random variables on unitary groups, the {I}tzykson-{Z}uber integral, and free probability, {\em Int. Math. Res. Not.} {\bf 17},953--982, (2003).

\bibitem{CD07} Collins B., Dykema K., A Linearization of Connes' Embedding Problem,
 {\em arXiv:0706.3918 } (2007).



\bibitem{cmss} Collins B., Mingo J., \'Sniady P. and Speicher R.,
Second Order Freeness and
Fluctuations of Random Matrices III.
Higher order freeness and free cumulants,
{\em Documenta Math.} (2007).

\bibitem{cs} Collins  B. and  Sniady P.,
Integration with respect to the Haar measure on unitary,
orthogonal and symplectic groups, {\tt math-ph/0402073}.
\bibitem{EM}{Ercolani, N. M. and McLaughlin K. D. T.-R.},
     {Asymptotics of the partition function for random matrices via
              {R}iemann-{H}ilbert techniques and applications to graphical
              enumeration},
   {\it  Int. Math. Res. Not.}{\bf 14} ,{755--820},({2003}).

\bibitem{FPEDZ}{Ferrer, A. , and Eynard, B., Di Francesco, P. and
              Zuber, J.-B.},
      {Correlation functions of {H}arish-{C}handra integrals over the
              orthogonal and the symplectic groups},
    {\em J. Stat. Phys.},
  {\bf 129}, {885--935},
      ({2007}).
  
\bibitem{G1}{Guionnet  A.},
      {First order asymptotics of matrix integrals; a rigorous
              approach towards the understanding of matrix models},
   {\it  Comm. Math. Phys.},
  {\bf 244}, {527--569},
      ({2004})
\bibitem{GM1} Guionnet A., Ma\"{\i}da  M.,
A {F}ourier view on the {$R$}-transform and related
              asymptotics of spherical integrals,
{\it J. Funct. Anal.}  {\bf 222},435--490,(2005).
\bibitem{GMa1}
     {Guionnet A. and  Maurel-Segala E.},
      {Combinatorial aspects of matrix models.},
   {\it Alea, electronic}
     ({2006}).
\bibitem{GMa2}
     {Guionnet A. and  Maurel-Segala E.},
      {Second order asymptotics for matrix models.},
   {\it   Ann. Probab.}  {\bf 35},2160--2212, (2007).

    
\bibitem{GZ} Guionnet A. and  Zeitouni O.,
      {Large deviations asymptotics for spherical integrals},
   {\it  J. Funct. Anal.},
 {\bf 188},{461--515} ,({2002}).


\bibitem{HC} Harish-Chandra,
      {Differential operators on a semisimple {L}ie algebra},
   {\it  Amer. J. Math.},
 {\bf 79},{87--120},({1957}).


\bibitem{Ho}  't Hooft G.,
 A planar diagram theory for strong interactions
  {\it   Nuclear Physics B}\textbf{72}461--473 (1974)
\bibitem{hp} F. Hiai and D. Petz, The semicircle law,
free random variables and entropy, AMS Mathematical Surveys and
Monographs 77, 2000.

   \bibitem{KS07} Klep I., Schweighofer M.,
   Connes' embedding conjecture and sums of hermitian squares
  {\it  Advances in Mathematics} \textbf{217}, No. 4, 1816-1837 (2008)


\bibitem{Ma}{ Maurel-Segala E.},
     {High order expansion for matrix models and enumeration
of maps .},  {\it  arviv },http://front.math.ucdavis.edu/0608.5192,
   ({2006})



\bibitem{HigherOrderFreeness2}{Mingo J., \'Sniady P. and Speicher R.
},
{Second order freeness and fluctuations of random matrices: {II.}},
  {U}nitary random matrices.
{\em Adv. in Math.} 209: 212 -- 240, ({2007}).

\bibitem{mst}{Mingo J., \'Sniady P., Speicher R. and Tan E.}
{Second Order Cumulants of products},
{\em preprint} arXiv:0708.0586, ({2007}).

\bibitem{NS}
  Nica, Alexandru and Speicher, Roland,
Lectures on the combinatorics of free probability,
{\it London Mathematical Society Lecture Note Series},
 {\bf 335},
(2006)


\bibitem{s1} Speicher R., Multiplicative functions on the lattice of
noncrossing partitions and free convolution, {\em Math. Ann.},
{\bf 298}, 611--628, (1994)
\bibitem{Zv}{Zvonkin A.},
     {Matrix integrals and map enumeration: an accessible
              introduction},
    {\it  Math. Comput. Modelling},
  {\bf 26},{281--304},({1997}).
\bibitem{Tu}
Tutte W. T, On the enumeration of planar maps, {\em
Bull. Amer. Math. Soc.},
{\bf  74}(1968), 64--74.
\bibitem{vo} Voiculescu  D.V., Lectures on free
probability theory, {\em Lecture Notes in Math.} {\bf 1738}, 279--349,
(2000).
\bibitem{vo2}
 Voiculescu D.,
A strengthened asymptotic freeness result for random matrices with
  applications to free entropy.
 {\em Internat. Math. Res. Notices}, (1):41--63, (1998).

\bibitem{Voi6} Voiculescu  D.V.,
The analogues of entropy and of {F}isher's information measure
              in free probability theory. {VI}. {L}iberation and mutual free
              information {\it Adv. Math.} {\bf 146}, 101--166,
(1999).
Lectures on free
probability theory, {\em Lecture Notes in Math.} {\bf 1738}
(2000), 279--349.
\bibitem{vdn} Voiculescu D.V.,  Dykema  K.J. and  Nica A.,
Free random variables, American Mathematical Society, Providence,
RI, (1992).
\bibitem{we} Weingarten D., Asymptotic behavior of group integrals in the limit of infinite rank, {\em J. Math. Phys.} {\bf 19}, 999--1001, (1978).
\bibitem{Z} Zuber J.B, On the large N limit of matrix integrals over the orthogonal grouparx, {\em arXiv:0805.0315 } (2008)
\bibitem{ZZ} Zinn Justin  P. and  Zuber J.B, On
some integrals over the $U(N)$ unitary group
and their large $N$ limit, {\em J. Phys. A}{\bf 36}, 3173--3194, (2003)
\end{thebibliography}
\end{document}